\documentclass[11pt]{article}
\usepackage{natbib}
\usepackage{bm}
\usepackage{amsthm,amsmath,amsfonts,amssymb}
\usepackage{tikz}
\usepackage[framemethod=latex]{mdframed}
\theoremstyle{Definition}

\newtheorem{Corollary}{Corollary}
\newtheorem{Example}{Example}
\newtheorem{Remark}{Remark}

\newtheorem{Proposition}{Proposition}
\newtheorem{Condition}{Condition}

\usepackage{geometry}
\usepackage{marginnote}

\def\gammab{\bar{\gamma}}
\def\varphib{\bar{\varphi}}

\def\alphab{\bar{\alpha}}
\def\thetab{\bar{\theta}}
\def\xb{\bar{x}}
\def\yb{\bar{y}}
\mdfdefinestyle{MyFrame}{%
linecolor=gray!15!white,
backgroundcolor=gray!40!white}

\mdfdefinestyle{MyFrame1}{%
 linecolor=gray!40!white,
 backgroundcolor=gray!9!white}

\mdfdefinestyle{MyFrame2}{%
   linecolor=gray!40!white,
  backgroundcolor=gray!0!white}

\title{\Large
A class of non-reversible hypercube long-range random walks and Bernoulli autoregression }
\date{\today}
 
\author{
Andrea Collevecchio\thanks{Andrea.Collevecchio@monash.edu, ORCID 0000-0001-6303-7925}\\School of Mathematical Sciences, Monash  University, Australia
\and
Robert C. Griffiths\thanks{Bob.Griffiths@monash.edu, ORCID 0000-0001-7190-5104, Corresponding Author}
\\School of Mathematical Sciences, Monash  University, Australia\\[0.1 cm]
}

\begin{document}
\maketitle
\begin{abstract}
We study a large class of  long-range random walks  which take values on the vertices
 of an $N$ dimensional hypercube.  These processes are  connected with multivariate Bernoulli autoregression. 
\end{abstract}
%{\bf Data Availability}\ The datasets generated during and/or analysed during the current study are available from the corresponding author on reasonable request.\\[0.1cm]
%

%%%%%%%%%%%%%%%%%%%%%%%%%%%%%
\section{Introduction}
\cite{CG2021} study a class of reversible long-range random walks $ (X_t)_{t}$  on the hypercube ${\cal V}_N =\{0,1\}^N$. The random walk is characterized by having an $N$-product Bernoulli$(p)$ stationary distribution and the property that 
\begin{equation*}
\mathbb{P}(X_t(B)\in C\mid X_{t-1}) = \mathbb{P}(X_t(B)\in C\mid X_{t-1}(B)),
%\label{subsets:00}
\end{equation*}
where $B\subseteq [N]=\{1,2,\ldots,N\}$, $C \subseteq \{0,1\}^B$ and $X_t(B)$ denotes a restriction to the coordinates in $B$.  The characterization, under Conditions \ref{condition:1} and \ref{condition:2} in Section \ref{Section:characterization}, is equivalent to a Lancaster characterization of the eigenvalues in the spectral expansion of the transition functions with tensor product eigenfunctions constructed from products of orthogonal functions on the marginal distributions. Knowing the spectral expansion allows a calculation of cut-off rates as $N \to \infty$. These mixing time rates are very fast in some cases, taking only two or three steps. In this paper we consider an extension of the processes in \cite{CG2021} to non-stationary long-range random walks. In the reversible model if $p_{xy}$, $x,y\in \{0,1\}$ is the marginal transition probability at a coordinate and $p_0=1-p, p_1=p$, then there is a reversibility property that $p_xp_{xy} = p_yp_{yx}$. Here this condition is relaxed somewhat in Condition \ref{condition:1} so that { \bf if} a coordinate $X_t[k]$ has a Bernoulli $(\gamma)$  distribution, {\bf then} $X_{t+1}[k]$ is distributed as Bernoulli$(\varphi)$, where $\gamma, \phi$ are specific parameters of the model. This extended concept of reversibility may seem unfamiliar, but commonly  holds  in $AR(1)$ series. Further comments about this assumption can be found in  Section \ref{Section:characterization}.
It is possible to construct $(X_t)_t$ with state space $\{0,1\}^\infty$ which has exchangeable coordinates such that for any $N$ coordinates $(X_t^{(N)})_t$ belongs to the class studied. This is an application of de Finetti's representation for countably infinite exchangeable random variables.

\subsection{Existing literature and structure of the paper} \cite{ES2020} study Bernoulli vector autoregressive models with a structure 
\[
X_t[k] = A_t[k]X_{t-1}[k]+B_t[k](1-X_{t-1}[k]),\ k \in [N],
\]
where $A_t, B_t$ are random vectors with entries in $\{0,1\}$, independent of $X_{t-1}$.
Our process $(X_t)_t$ falls into this class, with details in Corollary \ref{auto:00z}.
 It is natural to have multiple coordinates changing in a single time step in an autoregressive model, which corresponds to the long range structure of the random walk. The sample paths of the coordinates $(X_t[k] \colon k\in [N])$ are conditionally independent given the  sample paths of $(A_t,B_t)_t$.
Section \ref{Section:autoregression} expresses $(X_t)_t$ in this form, and finds the stationary means and covariances. 

Section \ref{Section:limits} finds the limiting distribution of $X_t$ as $t\to \infty$. In the limit when the process is time-homogeneous and $\varphi < \gamma$, $X_\infty \in {\cal V}_N$ is a Bernoulli vector with probabilities  
\begin{equation*}
\begin{aligned}
\mathbb{P}\big (X_\infty[k]=1\mid (S_l)_{l=0}^\infty) &= \varphi - (\gamma-\varphi)\sum_{l=1}^\infty\Big (\frac{\varphi}{\gamma}\Big )^l
\prod_{j=1}^l\big (1 - Z_j[k]/\varphi\big )\\
&=\varphi - (\gamma-\varphi)\sum_{l=1}^\infty\Big (\frac{\varphi}{\gamma}\Big )^l
\Big (-\frac{\varphib}{\varphi}\Big )^{S_l[k]},
\end{aligned}
\label{Intro:0}
\end{equation*}
with $(Z_j)_{j\in \mathbb{N}}$  \emph{i.i.d.} Bernoulli vectors where the joint distribution of  the coordinates of $Z_j$ on ${\cal V}_N$ are arbitrary and $S_l = \sum_{j=1}^lZ_j$. If $\varphi=\gamma$ then the distribution of $X_\infty[k]$ conditional on $(S_l)_{l=0}^\infty$, or unconditionally, is Bernoulli$(\varphi)$.
%%%%%%%
%
In Example \ref{deFinetti} the coordinates of $X_t$ are $N$ specific entries taken from an infinite exchangeable sequence $Y_t\in {\cal V}_\infty$.
% with a de Finetti measure $\nu_t$ such that for any $r\in \mathbb{N}$
%\begin{equation}
%\mathbb{P}\big (Y_t[1]=1,Y_t[2]=1,\ldots,Y_t[r]=1\big ) = \int_{[0,1]}\theta^r\nu_t(d\theta).
%\label{de:26}
%\end{equation}
%If $\nu_t=\nu$ is time homogeneous the limit distribution is shown  to be a Bernoulli vector $Y_\infty \in \{0,1\}^\infty$ where the conditional probability of coordinates being equal to one in $Y_\infty$ when $\varphi < \gamma$ and $\varphi+\gamma > 1$ is
%\begin{equation}
%\mathbb{P}\big (Y_\infty[k]=1 \mid (\theta_j)_{j=1}^\infty\big ) = \varphi - (\gamma-\varphi)\sum_{l=1}^\infty\big (\frac{\varphi}{\gamma}\big )^l\prod_{j=1}^l(1 - \theta_j/\varphi).
%\label{Intro:20}
%\end{equation}
%$(\theta_j)_{j=1}^\infty$ is an \emph{i.i.d.}  sequence of random variables in $[0,1]$  with measure $\nu$.
%(\ref{Intro:20}) has a similar form to (\ref{Intro:0}), with $Z_j[k]$ replaced by $\theta_j$, not depending on $k$. There is the distinction that $Z_j[k] \in \{0,1\}$ while $\theta_j \in [0,1]$.

Section \ref{Section:exchangeable} considers the case when $X_t$ has exchangeable coordinates. The transition distribution of $X_t$ is found to depend on the Hamming distance $\|x_t\|$, and also depends on $x_{t-1}$ through $\|x_{t-1}\|$ and  the inner product $\langle x_t,x_{t-1}\rangle$.
Proposition \ref{Hamming:000} shows that the Hamming distance $\|X_t\|$ is Markov and that
the transition density can be expressed as a diagonal expansion of Kratchouk
polynomials in $\|x_t\|$ and $\|x_{t-1}\|$. %The stationary distribution of the Hamming distance is found in Proposition \ref{Proposition:hamming:stationary} when the coordinates of $(Z_t)_t$ are exchangeable and $(Z_t)_t, (X_t)_t$ are time homogeneous.

In Section \ref{Section:normal} normal limit theorems for the transition density of $\|X^{(N)}_t\|$ in an exchangeable model are found. Proposition \ref{HammingNormal:01} shows that if $V^{(N)} = N^{-1/2}\big (\|X^{(N)}_t\| - N\varphi\big )$ and $N^{1/2}(\varphi^{(N)} - \gamma^{(N)}) \to c$, a constant, then the limit process $\{V_t\}$ is a (normal) AR(1) series with random parameters.
If the difference between the parameters is of constant order instead, then $V^{(N)} \to \pm \infty$, showing that it is appropriate to model the difference as small for $N$ large to obtain stability.

In Section \ref{Section:estimation} a sketch is made of how to estimate the distribution of $\varphi,\gamma$ and the distribution of $Z_t$ in a time-homogeneous exchangeable model.

\begin{Remark}[Notation]
We denote $\bar{a} = 1- a$, for any $a \in (0,1)$, set $\alpha := \min(\varphi, \gamma)$, and $\psi := \varphi/\gamma$ if $\varphi \leq \gamma$ or $\psi := \varphib/\gammab$ if $\varphi > \gamma$.
%\[
%\psi := %\frac{\alpha\varphib}{\alphab\gamma}=
%\begin{cases}
%\frac{\varphi}{\gamma},&\varphi \leq \gamma\\
%\frac{\varphib}{\gammab},&\varphi > \gamma 
%\end{cases}.
%\]
\end{Remark}
\section{Characterization and spectral expansion}\label{Section:characterization}
We first focus on the transition of a  single coordinate  from  $x$ to $y$,  where $x,y\in \{0,1\}$.  It is assumed  that when  $X$ has a Bernoulli($\gamma$) distribution then  $Y$ has a Bernoulli($\varphi$) distribution. This transition can be expressed when $\varphi_1=\varphi$, $\varphi_0=1-\varphi$, as
\begin{equation}
p_{xy} = \varphi_y\{1 + \kappa (1-x/\gamma)(1-y/\varphi)\},
\label{little:00}
\end{equation}
where $\kappa$ is a constant such that $p_{xy}\geq 0$. 
  This is a spectral expansion in one transition using the  orthogonal function set on $Y$ of $\{1,1-y/\varphi\}$ and similarly for $X$. The joint distribution of $(X,Y)$ always has an orthogonal function expansion
\[
\phi_y\gamma_x\big \{1 +  \kappa_{01}(1-y/\varphi) +  \kappa_{10}(1-x/\gamma)+
\kappa (1-x/\gamma)(1-y/\varphi)\big \}
\]
for  $\gamma_1 = \gamma = 1- \gamma_0$ and constants $\kappa_{01},\kappa_{10},\kappa$. Insisting that the marginal distribution of $Y$ is Bernoulli$(\varphi)$ forces $\kappa_{01}=0$ and similarly $\kappa_{10}=0$ since $X$ is Bernoulli$(\gamma)$.

Let $\alpha = \min \{\varphi, \gamma\}$. It is straightforward to show that $p_{xy}\geq 0$ is equivalent to $-1 \leq \kappa \leq \alpha/\alphab$ if $\varphi+\gamma \geq 1$, and   $-\varphi\gamma/\varphib\gammab \leq \kappa \leq \alpha/\alphab$ if $\varphi + \gamma < 1$.  Notice that if $\varphi+\gamma < 1$ then $1-\varphi + 1 - \gamma > 1$, so if appropriate $x,y$ might be exchanged for $\xb,\yb$ and $\varphi,\gamma$ for $\varphib,\gammab$. In the paper we take $\varphi+\gamma \geq 1$; if $\varphi=\gamma$ this is equivalent to $\varphi \geq 1/2$.\\
%
 
  %The $2\times 2$ transition matrix 
%  $$
% \begin{pmatrix}
%p_{00} & p_{01}\\
%p_{10} & p_{11}
%\end{pmatrix}
%$$
  Any non-degenerate markov chain on $\{0, 1\}$ is reversible and ergodic.
 In our case, where the transition probabilities are determined by \eqref{little:00},  the stationary measure is   Bernoulli with parameter
\begin{equation}
\frac{1-\kappa(1-\varphi)}{1 - (1-\varphi - \gamma)\kappa/\gamma}.
\label{little:01}
\end{equation}
%\marginnote{\textcolor{purple}{Iff condition for the hypercube to be reversible?}}
%The assumption (\ref{little:00}) is additional to having an arbitary $2$-dimensional Markov chain for given $\varphi,\gamma$.
 Coordinates of the hypercube that we study are generally not independent, the usual reversibility property does not hold jointly for the coordinates, although it does hold marginally for the coordinates, and the stationary distribution is not a product distribution unless $\varphi=\gamma$. The assumption (\ref{little:00}) does imply that the marginal stationary distributions are identical Bernoulli with parameter (\ref{little:01}).
 In the spectral expansion of the transition density of the long-range hypercube random walk we consider, the eigenfunctions belong to tensor product sets of the orthogonal functions at coordinates. That is, functions in 
$\otimes_{k=1}^N\{1,1-y[k]/\varphi\}$, and similarly for $x$.
$(X_t)_t$ will be said to belong to a class $\mathcal{G}$ if it satisfies the following two Conditions \ref{condition:1}, \ref{condition:2}.
\begin{mdframed}
[style=MyFrame1]
\begin{Condition}\label{condition:1} \normalfont
For each $t\in \mathbb{N}$, $X_t$ has a product-measure
$\otimes_{i=1}^N \Phi_i$ when $X_{t-1}$ has a product-measure $\otimes_{i=1}^N \Gamma_i$, where $\Phi_i$ and $\Gamma_i$ are respectively Bernoulli$(\varphi)$ and Bernoulli$(\gamma)$ measures.
We can assume, without loss of generality that $\varphi + \gamma \geq 1$.  If the latter does not hold, we can achieve it by reversing the roles  of $1$ and $0$.  We assume the following  form of reversibility
\begin{equation}
\begin{aligned}
\gamma^{\|\bm{x}\|}\gammab^{N-\|\bm{x}\|}&\mathbb{P}\big (X_t=\bm{y}\mid X_{t-1}=\bm{x};\varphi,\gamma\big )= \varphi^{\|\bm{y}\|}\varphib^{N-\|\bm{y}\|}
 \mathbb{P}\big (X_t=\bm{x}\mid X_{t-1}=\bm{y};\gamma,\varphi\big ).
\end{aligned}
\label{psuedoreverse:00}
\end{equation}
Notice how the parameters $\varphi,\gamma$ are interchanged across the sides of (\ref{psuedoreverse:00}).
If $\varphi=\gamma$ then (\ref{psuedoreverse:00}) coincides with the canonical reversibility condition. 
\end{Condition}
\end{mdframed}
We can rewrite (\ref{psuedoreverse:00}) as follows
\begin{equation}
\mathbb{E}\Big [f(\bm{Y})g(\bm{X});\varphi, \gamma)\Big ]
= \mathbb{E}\Big [f(\bm{X})g(\bm{Y});\gamma,\varphi)\Big ].
\label{fgeq:00}
\end{equation}
when the marginal distributions of $\bm{Y},\bm{X}$ are product Bernoulli $\varphi,\gamma$.
%\marginnote{\textcolor{purple}{Mass transfer concept?}}

\begin{Remark}
 AR(1) models satifsy a condition similar Condition to \ref{condition:1}.
Let $X_t=\mu + \lambda X_{t-1} + \epsilon_t$, where $X_{t-1}$ is N$(a_0,\kappa_0^2)$ and $\epsilon_t$ is independent N$(0,\sigma^2)$. Then $X_t$ is N$(a_1,\kappa_1^2)$ where
$a_1=\mu+\lambda a_0$, $\kappa_1^2 = \lambda^2\kappa_0+\sigma^2$. The covariance between $X_t$ and $X_{t-1}$ equals to $v=\lambda\kappa_0^2$. Let $f(y\mid x; a_1,a_0,\kappa_1,\kappa_0,v)$
be the conditional density of $X_t$ given $X_{t-1}=x$ with an ordered parameter list.
Then (\ref{psuedoreverse:00}) is analogous to 
\begin{equation}
\begin{aligned}
&\frac{1}{\sqrt{2\pi}\kappa_0}e^{-\frac{1}{2}(x-a_0)^2/\kappa_0^2}
f(y\mid x; a_1,a_0,\kappa_1,\kappa_0,v) = \frac{1}{\sqrt{2\pi}\kappa_1}e^{-\frac{1}{2}(y-a_1)^2/\kappa_1^2}f(x\mid y;a_0,a_1,\kappa_0,\kappa_1,v),
\label{ARRev:00}
\end{aligned}
\end{equation}
which follows from the form of the bivariate normal distribution.
Notice that again the order of the parameter list is changed from the left side to the right side.
\end{Remark}
Generally in stochastic processes reversibility provides a description of  stationarity  of the process. That is not the case here,  (\ref{psuedoreverse:00}) and (\ref{ARRev:00}) are just statements about what happens in one transition.
\begin{mdframed}
[style=MyFrame1]
\begin{Condition}\label{condition:2}\normalfont For any  $\bm{y} \in {\cal V}_N$, for all $B\subseteq [N]$, let $\bm{y}(B)$ be the  projection from ${\cal V}_N$ on $B$ defined as  the vector $\bm{y}(B) = (\bm{y}[j], j \in B)$. For each $t\in \mathbb{N}$ we assume a local change condition that for all $C \subseteq \{0,1\}^B$
\begin{equation}
\mathbb{P}(X_t(B)\in C\mid X_{t-1}) = \mathbb{P}(X_t(B)\in C\mid X_{t-1}(B)).
\label{subsets:0}
\end{equation}
\end{Condition}
\end{mdframed}
 A characterization of the spectral expansion for the hypercube transition density is found when $(X_t)_t\in \mathcal{G}$ in Proposition \ref{Proposition:03}. The following is an extension of Proposition 1 of \citet{CG2021}, when $\varphi = \gamma$.
\begin{mdframed}
[style=MyFrame1]
\begin{Proposition}\label{Proposition:03} The process
$\bm{X} \in {\mathcal G}$ if and only if there exists a random variable $Z_t$  in ${\cal V}_N$ such that 
$p_{t-1\to t}(\bm{x},\bm{y}):= \mathbb{P}\big (X_t=\bm{y}\mid X_{t-1}=\bm{x}\big )$ can be expanded as
\begin{equation}
p_{t-1\to t}(\bm{x},\bm{y}) = \varphi^{\|\bm{y}\|}\varphib^{N-\|\bm{y}\|}\Bigg \{
1 + \sum_{A\subseteq [N], A\ne \emptyset}\kappa_{t,A}\prod_{k \in A}
\Big (1 - \frac{\bm{y}[k]}{\varphi}\Big )\Big (1 - \frac{\bm{x}[k]}{\gamma}\Big )\Bigg \},
\label{hypercube:000}
\end{equation}
where
\begin{equation}
\begin{aligned}
\kappa_{t,A}&= \Big (\frac{\alpha}{\alphab}\Big )^{|A|}\mathbb{E}\Big [
\prod_{k \in A}\Big (1 - \frac{Z_t[k]}{\alpha}\Big )
\Big ]\\
&= 
(-1)^{|A|}\mathbb{E}\Big [\Big (-\frac{\alpha}{\alphab}\Big )^{\sum_{k\in A}(1-Z_t[k])}\Big ],
\end{aligned}
\label{kappa:00}
\end{equation}
with $\alpha=\min\{\varphi,\gamma\}$.
\end{Proposition}
\end{mdframed}
\begin{proof}
%The proof is an extended version of the proof in \cite{CG2021} when $\varphi=\gamma$.\\[0.1cm]
\emph{Sufficiency.}
Suppose (\ref{hypercube:000}) holds. \\
\emph{Condition \ref{condition:1}}.
The product-measure requirement is satisfied because, under the product-measure $\otimes_{i=1}^N \Gamma_i$,
\[
\mathbb{E}\Big [
\prod_{k \in A}
\Big (1 - \frac{\bm{X}[k]}{\gamma}\Big )\Big ]=0.
\]
Condition \ref{condition:1}  is easily checked to be true, taking care with the parameter reversal.\\
\emph{Condition \ref{condition:2}}.
The marginal distribution of $\bm{Y}(B)\mid X_{t-1}$ is seen to be 
\begin{equation}
\varphi^{\|\bm{y}(B)\|}\varphib^{N-\|\bm{y}(B)\|}\Bigg \{
1 + \sum_{A\subseteq [B], A\ne \emptyset}\kappa_{t,A}\prod_{k \in A}
\Big (1 - \frac{\bm{y}[k]}{\varphi}\Big )\Big (1 - \frac{\bm{x}[k]}{\gamma}\Big )\Bigg \},
\label{hypercube:000B}
\end{equation}
noting that in a product-measure $\otimes_{i\in [N]-B}\Phi_i$, 
\[
\mathbb{E}\Big [\prod_{k \in A}
\Big (1 - \frac{\bm{y}[k]}{\varphi}\Big )\Big ]=0\text{~if~}A \not\subset B.
\]
Condition \ref{condition:2} is satisfied because (\ref{hypercube:000B}) only depends on $\bm{x}(B)$.\\[0.1cm] 
\emph{Necessity.}
Suppose Condition \ref{condition:1} and Condition \ref{condition:2} are satisfied.\\
 Set $\bm{X} = X_{t-1}$, $\bm{Y} = X_t$. Since $\bm{X},\bm{Y}\in {\cal V}_N$ there are constants $\{\kappa_{BA}\}_{  B,A\subseteq [N]}$ such that 
\[
\mathbb{E}\Big [
\prod_{k\in B}\Big ( 1 - \frac{\bm{Y}[k]}{\varphi}\Big )\>\Big |\> \bm{X}\Big ]
= \sum_{A\subseteq [N]}\Big (\frac{\varphib}{\varphi}\Big )^{|A|}\kappa_{BA}\prod_{l\in A}\Big ( 1 - \frac{\bm{X}[l]}{\gamma}\Big ).
\]
Let $\bm{X}$ have a product-measure $\otimes_{i=1}^N \Gamma_i$.
If $B\ne \emptyset$ then $\kappa_{B\emptyset}=0$ because of the product-measure criterion in Condition \ref{condition:1}. Also
\[
\mathbb{E}\Big [
\prod_{k\in B}\Big ( 1 - \frac{\bm{Y}[k]}{\varphi}\Big )
\prod_{l\in A}\Big ( 1 - \frac{\bm{X}[l]}{\gamma}\Big ) \Big ]
= \kappa_{BA}\Big (\frac{\varphib}{\varphi}\Big )^{|B|}\Big (\frac{\gammab}{\gamma}\Big )^{|A|}.
\]
$\kappa_{BA} \ne 0$ if $A\not\subset B$ is inconsistent with (\ref{subsets:0}), therefore $\kappa_{BA}=0$ if $A \not\subset B$. Let $f,g$ be bounded functions from ${\cal V}_N$ to $\mathbb{R}$.
Taking 
\[
f(\bm{Y}) = \prod_{k\in B}\Big ( 1 - \frac{\bm{Y}[k]}{\varphi}\Big ),\>
g(\bm{X}) = \prod_{l\in A}\Big ( 1 - \frac{\bm{X}[l]}{\gamma}\Big )
\]
in (\ref{fgeq:00}) shows that 
\[
\kappa_{BA}\Big (\frac{\varphib}{\varphi}\Big )^{|B|}\Big (\frac{\gammab}{\gamma}\Big )^{|A|}
=\kappa_{AB}\Big (\frac{\varphib}{\varphi}\Big )^{|A|}\Big (\frac{\gammab}{\gamma}\Big )^{|B|}
\]
and therefore $\kappa_{BA} = 0$ if $A\ne B$. Set $\kappa_{t,A} := \kappa_{AA}$.
It will now be shown that (\ref{kappa:00}) holds. Directly from (\ref{hypercube:000})
\[
\mathbb{E}\Big [\prod_{k \in A}\Big (1 - \frac{\bm{Y}[k]}{\varphi}\Big )\>\Big |\> \bm{x}\Big ]
= \Big (\frac{\varphib}{\varphi}\Big )^{|A|}\kappa_{t,A}
\prod_{k \in A}\Big (1 - \frac{\bm{x}[k]}{\varphi}\Big ).
\]
Therefore it is necessary that $\exists$ $Z_t \in \{0,1\}^N$, by setting $Z_t = \bm{Y}\mid \bm{x}=\bm{0}$, such that
\begin{equation}
\kappa_{t,A}=
\Big (\frac{\varphi}{\varphib}\Big )^{|A|}\mathbb{E}\Big [
\prod_{k \in A}\Big (1 - \frac{Z_t[k]}{\varphi}\Big )
\Big ].
\label{mixture:a}
\end{equation}
Similarly it is also necessary that $\exists$ $Z^\prime_t \in \{0,1\}^N$ such that
\begin{equation}
\kappa_{t,A}=
\Big (\frac{\gamma}{\gammab}\Big )^{|A|}\mathbb{E}\Big [
\prod_{k \in A}\Big (1 - \frac{Z_t^\prime[k]}{\gamma}\Big )
\Big ].
\label{mixture:b}
\end{equation}
by taking $Z_t^\prime = \bm{X}\mid \bm{y}=0$, when $\bm{X}$ has product-measure $\otimes_{i=1}^N \Gamma_i$.
This completes the proof of the necessity by choosing (\ref{mixture:a}) or (\ref{mixture:b}) according to whether $\varphi \leq \gamma$ or $\varphi>\gamma$.
\end{proof}
%%$
\begin{mdframed}
[style=MyFrame1]
\begin{Corollary}\label{Corr:pt}
Recall that $p_{t-1\to t}(\bm{x},\bm{y}):= \mathbb{P}\big (X_t=\bm{y}\mid X_{t-1}=\bm{x}\big )$. We have
\begin{equation}
\begin{aligned}
p_{t-1\to t}(\bm{x},\bm{y}) &= \varphi^{\|y\|}\varphib^{N-\|y\|}
\mathbb{E}\Bigg [\prod_{k=1}^N
\Big \{ 1 + \frac{\alpha}{\alphab}\Big ( 1 - \frac{Z_t[k]}{\alpha}\Big )\Big ( 1 -\frac{\bm{y}[k]}{\varphi}\Big )\Big (1 -\frac{\bm{x}[k]}{\gamma}\Big )\Big \}
\Bigg]\\
& =\varphi^{\|y\|}\varphib^{N-\|y\|}
\mathbb{E}\Bigg [\prod_{k=1}^N
\Big \{ 1 - \big (-\frac{\alphab}{\alpha}\big )^{1-Z_t[k]}
\big (-\frac{\varphib}{\varphi}\big )^{\bm{y}[k]}
\big (-\frac{\gammab}{\gamma}\big )^{\bm{x}[k]}
\Big \}\Bigg],
\end{aligned}
\label{hypercube:001}
\end{equation}
with notation as in Proposition \ref{Proposition:03}.
\end{Corollary}
\end{mdframed}
\begin{proof}
Expand the product in (\ref{hypercube:001}) to see the equivalence with (\ref{hypercube:000}).
\end{proof}
%\marginnote{\textcolor{purple}{Maybe cross reference with (\ref{little:00})}}
Notice that $\kappa_{t,A}=0$ for all $A$ implies that  $(X_t)_t$ are  independent random variables. It is  also a necessary condition and is equivalent to $Z_t$ having all coordinates \emph{i.i.d.} Bernoulli$(\alpha)$. Then,  for each $t$, $X_t$ has \emph{i.i.d.} Bernoulli$(\varphi)$ coordinates. If $Z_t$ has $N-r$ \emph{i.i.d.} Bernoulli$(\alpha)$ coordinates then  $\kappa_{t,A}=0$ if $|A|>r$.

More generally, $(X_t)_t$ follows a hidden Markov model, controlled by $(Z_t)_t$. 
The coordinates of $(X_t)_{t}$, conditionally on $(Z_t)_{t}$, are independent. $Z_t$ can be chosen to have any multivariate Bernoulli distribution, for example those developed in \citet{T1990,FS2018}.
We do emphasise though that the existence of $(Z_t)_t$ is a \emph{consequence} of Conditions \ref{condition:1} and  \ref{condition:2}. 
In both (\ref{hypercube:000}) and (\ref{hypercube:001}) the conditional distribution of $\bm{Y}$ given $\bm{X}$ and $Z_t$ is similar, with $\mathbb{E}$ with respect to $Z_t$ removed.
\begin{mdframed}
[style=MyFrame1]
\begin{Corollary}
The Lancaster problem of characterizing which eigenvalue sequences $\{\kappa_{t,A}\}_{A\subseteq [N]}$ make (\ref{hypercube:000}) non-negative and therefore a probability distribution is solved by (\ref{kappa:00}).
\end{Corollary}
\end{mdframed}
\begin{Remark}
$\mathcal{G}$ is a convex set. Let $p^\circ_{t-1\to t}(\bm{x},\bm{y})$ and $p^\prime_{t-1\to t}(\bm{x},\bm{y})$ be two transition densities satisfying (\ref{hypercube:000}) with (\ref{kappa:00}) holding for two random variables $Z_t^\circ$ and $Z_t^\prime$. For $\lambda \in [0,1]$ let
$Z_t = Z_t^\circ$ with probability $\lambda$ and $Z_t = Z_t^\prime$ with probability $1-\lambda$. Then the transition density corresponding to $Z_t$ is  $\lambda p^\circ_{t-1\to t}(\bm{x},\bm{y}) +  (1-\lambda)p^\prime_{t-1\to t}(\bm{x},\bm{y})$, showing convexity. Alternatively convexity follows directly from Conditions \ref{condition:1} and \ref{condition:2}.
The extreme points of $\mathcal{G}$ correspond to non-random $z_t \in {\cal V}_N$ and any $(X_t)_t\in \mathcal{G}$ has a mixture distribution of the extreme points. Extreme points of $(X_t)_t \in \mathcal{G}$ with exchangeable coordinates correspond to $Z_t$ with a non-random $M_t \in [N]$ entries one, exchangeably distributed. 
\end{Remark}
Let $(U(t, k))_{t, k}$ be a collection of\emph{ i.i.d.} uniform random variables. 
\begin{mdframed}
[style=MyFrame1]
\begin{Proposition}\label{construction:250}
$\bm{X} \in {\mathcal G}$ if and only if it can be constructed in the following way. Let $Z_t \in \{0,1\}^N$ and take $\varphi + \gamma \geq 1$, without loss of generality. In a transition from $t-1$ to $t$ independently change coordinates $k \in [N]$ by the rules that
\begin{itemize}
\item[(a)] If $Z_t[k] =1$ then
\[
X_t[k]=
\begin{cases}
1&\text{if either ~}X_{t-1}[k]=0  \text{ or } ( X_{t-1}[k] =1 \mbox{ and } U(t, k) \ge  \varphib/\gamma)\\
0&\text{otherwise}.
\end{cases}
\]
\item[(b)] 
\begin{itemize} 
\item[(i)] If $Z_t[k] =0$ and $\varphi < \gamma$
  then
\[
X_t[k] =
\begin{cases}
0&\text{if~}X_{t-1}[k]=0 \text{ or } ( X_{t-1}[k] =1 \mbox{ and } U(t, k) \ge  \psi)\\
1&\text{otherwise}.
\end{cases}
\]
\item[(ii)] If $Z_t[k] =0$ and $\varphi > \gamma$
 then
\[
X_t[k]= 
\begin{cases}
1&\text{if~}X_{t-1}[k]=1 \text{ or } ( X_{t-1}[k] =0 \mbox{ and } U(t, k) \ge  \psi)\\
0&\text{otherwise}.
\end{cases}
\]
\end{itemize}
\end{itemize}
\end{Proposition}
\end{mdframed}
\begin{proof} 
Denote $Y=X_t$, $X=X_{t-1}$ in (a) and (b). 
Expansions of the transition probabilities at coordinate $k$ according to whether $Z_t[k]=1$ or $Z_t[k]=0$ are
\begin{eqnarray}
p_{1;\bm{x}[k]\>\bm{y}[k]}&=&\varphi^{\bm{y}[k]}\varphib^{1-\bm{y}[k]}\Big \{ 1 - \Big (1 -\frac{\bm{y}[k]}{\varphi}\Big )\Big (1 -\frac{\bm{x}[k]}{\gamma}\Big )\Big \}\nonumber \\
p_{0;\bm{x}[k]\>\bm{y}[k]} &= &\varphi^{\bm{y}[k]}\varphib^{1-\bm{y}[k]}
\Big \{ 1 + \frac{\alpha}{\alphab}\Big ( 1 -\frac{\bm{y}[k]}{\varphi}\Big )\Big (1 -\frac{\bm{x}[k]}{\gamma}\Big )\Big \}.
\label{singletrans:00}
\end{eqnarray}
The transition probability given $Z_t = \bm{z}$ is then
\begin{equation}
\begin{aligned}
\bm{z}[k]p_{1;\bm{x}[k]\>\bm{y}[k]} &+ (1-\bm{z}[k])p_{0;\bm{x}[k]\>\bm{y}[k]}\\
&= \varphi^{\bm{y}[k]}\varphib^{1-\bm{y}[k]}
\Big \{ 1 + \frac{\alpha}{\alphab}\Big ( 1 - \frac{\bm{z}[k]}{\alpha}\Big )\Big ( 1 -\frac{\bm{y}[k]}{\varphi}\Big )\Big (1 -\frac{\bm{x}[k]}{\gamma}\Big )\Big \}.
\end{aligned}
\label{proofz:00}
\end{equation}
To show that (a) and (b) are equivalent to (\ref{singletrans:00}) it is enough to verify that
 $\mathbb{E}\big [Y[k]\mid \bm{x}[k],\bm{z}[k]\big ]$ agree.  Suppose $Z[k]=1$. From (a)
 \begin{equation}
 \begin{aligned}
 \mathbb{E}\big [Y[k]\mid \bm{x}[k],Z[k]=1\big ] &=
 1 - \bm{x}[k] + (1 - \varphib/\gamma)\bm{x}[k]
 \end{aligned}
 \label{a:00000}
 \end{equation}
 and from  (\ref{singletrans:00})
 \[
 \begin{aligned}
  \mathbb{E}\big [Y[k]\mid x[k],Z[k]=1\big ] &=
  \varphi - \varphi(1-1/\varphi)(1-\bm{x}[k]/\gamma)\\
  &= 1 - (\varphib/\gamma)\bm{x}[k],
 \end{aligned}
 \]
 in agreement with (\ref{a:00000}). A similar argument shows agreement when $Z[k]=0$.

Since transitions at coordinates are made independently conditional on $Z_t=\bm{z}$ (\ref{proofz:00}) implies (\ref{hypercube:001}) which is equivalent to (\ref{hypercube:000}).
% The distribution of $Z_t$ is uniquely determined by (\ref{hypercube:000}) as that of $Y\mid \bm{x}=\bm{0}$.
\end{proof}

Proposition \ref{construction:250} shows that $p_{t-1\to t}(\bm{x},\bm{y})$ can be expressed as a mixture of independent chains at each coordinate  in $[N]$. For fixed $k$, the process  $(X_t[k])_t$ is a (non-homogeneous) Markov Chain with  transition probability matrix
\begin{equation*}
P_t =
\begin{pmatrix}
\varphib + \frac{\alpha\varphib}{\alphab}\Big ( 1 -\frac{1}{\alpha}\mathbb{E}\big [Z_t[k]\big ]\Big )
&\varphi - \frac{\varphib\alpha}{\alphab}
\Big ( 1 -\frac{1}{\alpha}\mathbb{E}\big [Z_t[k]\big ]\Big )\\
\varphib - \frac{\varphib\alpha\gammab}{\alphab\gamma}\Big ( 1 -\frac{1}{\alpha}\mathbb{E}\big [Z_t[k]\big ]\Big )
&
\varphi + \frac{\varphib\alpha\gammab}{\alphab\gamma}\Big ( 1 -\frac{1}{\alpha}\mathbb{E}\big [Z_t[k]\big ]\Big )
\end{pmatrix}.
\label{Ptk:00}
\end{equation*}
%\marginnote{\textcolor{purple}{Check this matrix and below.}}
If $(Z_t[k])_{t}$ is time-homogeneous  then $(X_t[k])_{t}$ is reversible with a Bernoulli $(p)$ stationary distribution,
  where
\begin{equation}
\begin{aligned}
p &= \frac
{
\varphi - \frac{\varphib\alpha}{\alphab}
\Big ( 1 -\frac{1}{\alpha}\mathbb{E}\big [Z_t[k]\big ]\Big )
}
{
1 - \frac{\alpha\varphib}{\alphab\gamma}\Big ( 1 -\frac{1}{\alpha}\mathbb{E}\big [Z_t[k]\big ]\Big )
}\ .
%\\
%&=\begin{cases}
%\frac{\mathbb{E}\big [Z_t[k]\big ]}
%{1 - \frac{\varphi}{\gamma}+\frac{1}{\gamma}\mathbb{E}\big [Z_t[k]\big ]},&\alpha = \varphi\\
%1 - \frac{\frac{\varphib}{\gamma} \mathbb{E}\big [Z_t[k]\big ]}
%{1 - \frac{\varphib}{\gammab}+ \frac{\varphib}{\gammab\gamma}\mathbb{E}\big [Z_t[k]\big ]},
%&\alpha=\gamma.
%\end{cases}
\end{aligned}
\label{Bernp:00}
\end{equation}
Although the single coordinate chains are reversible the full chain $(X_t)_{t}$ will not be reversible if $\varphi\ne \gamma$ and the coordinates of $(Z_t)_{t}$ are not independent. A simple case when the coordinates of $X_t$ \emph{are} independent is when $Z_t$ are \emph{i.i.d.} Bernoulli $(\omega)$ for each $t$.

Proposition \ref{Proposition:03} characterizes eigenvalues of $p_{t-1\to t}(\bm{x},\bm{y})$ in terms of $Z_t$ and Proposition \ref{construction:250} allows a probabilistic construction of $(X_t)_t$ from $(Z_t)_t$. The next Corollary provides a direct calculation of $\mathbb{P}(Z_t=\bm{z})$ in terms of $p_{t-1\to t}(\bm{z},\bm{0})$ or $p_{t-1\to t}(\bm{0},\bm{z})$.
\begin{mdframed}
[style=MyFrame1]
\begin{Corollary}
Let $(X_t)_t\in {\mathcal G}$. The distribution of $Z_t$ satisfies 
\begin{equation}
\mathbb{P}\big (Z_t = \bm{z}\big ) =
\begin{cases}
p_{t-1\to t}(\bm{z},\bm{0})&\text{if~}\varphi \leq \gamma\\
\Big (\frac{\gammab}{\varphib}\Big )^Np_{t-1\to t}(\bm{0},\bm{z})
&\text{if~}\varphi > \gamma
\end{cases}.
\label{construction:510}
\end{equation}
\end{Corollary}
\end{mdframed}
\begin{proof}
The distribution of $Z_t$ always has an orthogonal function expansion
\begin{equation}
\mathbb{P}\big (Z_t = \bm{z}\big )  = \alpha^{\|\bm{z}\|}\alphab^{N-\|\bm{z}\|}\Bigg \{
1 + \sum_{A\subseteq [N], A\ne \emptyset}\kappa_{t,A}
\prod_{k \in A}
\Big (1 - \frac{\bm{z}[k]}{\alpha}\Big )\Bigg \},
\label{hypercube:359}
\end{equation}
where $\kappa_{t,A}$ is given by (\ref{kappa:00}). Comparing (\ref{hypercube:359}) with (\ref{hypercube:000}) proves (\ref{construction:510}).
\end{proof}
\section{Bernoulli autoregressive system}\label{Section:autoregression}
\begin{mdframed}
[style=MyFrame1]
\begin{Corollary}\label{auto:00z}
$(X_t)_{t}$ satisfies  a Bernoulli autoregressive system
\begin{equation*}
X_{t}[k] = A_t[k]X_{t-1}[k]+B_t[k](1-X_{t-1}[k]),\quad k\in [N],
%\label{AR1:11}
\end{equation*}
where $A_t, B_t$ depend on $Z_t$ such that
\[
(A_t[k],B_t[k]) =
\begin{cases}
\begin{cases}
(I_{U(t, k) \le \varphi/\gamma} ,0)& \mbox{ if } \varphi \leq \gamma\\
(1,I_{U(t, k) > \varphib/\gamma})&\mbox{ if } \varphi > \gamma
\end{cases}
&\text{in the case~}Z_{t}[k]=0,\\
(I_{U(t, k) > \varphib/\gamma},1)&\text{in the case ~}Z_{t}[k]=1.
\end{cases}
\]

%
%If $\varphi=\gamma$, then 
%\[
%(A_t[k],B_t[k]) =
%\begin{cases}
%(1,0) &\text{if~}Z_{t}[k]=0,\\
%(I_{U(t, k) > \varphib/\varphi},1)&\text{if~}Z_{t}[k]=1.
%\end{cases}
%\]
%whith $\{\vartheta_{t}\}$ Bernoulli $(2-1/\varphi)$.
\end{Corollary}
\end{mdframed}
%The behaviour of how $X_t$ changes with $Z_t$ is different when $\varphi=\gamma$ and $\varphi\ne \gamma$. If $\varphi=\gamma$ then $Z_t[k]=0$ implies that $X_t[k]=X_{t-1}[k]$, however if $Z_t[k]=1$ then 
%$X_t[k] = I_{U(t, k) > \varphib/\gamma} X_{t-1}[k] + 1-X_{t-1}[k]$. If $\varphi \ne \gamma$ and $Z_t[k]=0$ then $X_t[k]=X_{t-1}[k]$ with probability $\varphi/\gamma$.

 It is straightforward to show that
\[
\text{Cov}\big (A_t[k],B_t[k]\big ) = 
\begin{cases}
-(\gammab/\gamma)
\text{Var}\big (Z_{t}[k]\big )&\varphi \leq \gamma\\
-(\varphib^2/(\gamma\gammab))
\text{Var}\big (Z_{t}[k]\big )&\varphi > \gamma.
\end{cases}
\]
$A_t$ and $B_t$ are independent if and only if $\max \{\varphi,\gamma\} =1$, when the coefficients of $\text{Var}\big (Z_{t}[k]\big )$  are zero, or $Z_{t}$ is a constant in ${\cal V}_N$ with probability 1.

The conditional means and variances in a transition are now calculated. From (\ref{hypercube:001})
\begin{equation}
\begin{aligned}
&\mathbb{E}\big [X_t[k]-\varphi \mid X_{t-1},Z_t\big ]
= -\varphib/(\alphab\gamma)(Z_t[k]-\alpha)(X_{t-1}[k]-\gamma),\\
&\mathbb{E}\big [X_t[k]\mid X_{t-1}\big ]
= (\varphi-\alpha)/\alphab +(\varphib/\alphab)\mathbb{E}\big [Z_t[k]\big ]+ (\alpha\varphib)/(\alphab\varphi)X_{t-1}[k] 
- \varphib/(\alphab\gamma)\mathbb{E}\big [Z_t[k]\big ]X_{t-1}[k],\\
&\text{Var}\big (X_t[k]\mid X_{t-1}\big )
= \mathbb{E}\big [X_t[k]\mid X_{t-1}\big ]\Big (1 - \mathbb{E}\big [X_t[k]\mid X_{t-1}\big ]\Big ).
\end{aligned}
\label{mean:000}
\end{equation}
For $k\ne l$
\begin{equation*}
\begin{aligned}
&\text{Cov}\big (X_t[k],X_t[l]\mid X_{t-1})\\
&~~=
0+\text{Cov}\big (\mathbb{E}\big [X_t[k]\mid X_{t-1},Z_t\big ],\mathbb{E}\big [X_t[l]\mid X_{t-1},Z_t\big ]\big )\\
&~~= \varphib/(\alphab\gamma)^2\text{Cov}\big (Z_t[k],Z_t[l]\big )
(X_{t-1}[k]-\gamma)(X_{t-1}[l]-\gamma)
\end{aligned}
\end{equation*}
In a process where $Z_t$ is time homogeneous and $(X_t)_{t}$ is stationary denote
\[
\begin{aligned}
&m_Z=\mathbb{E}\big [Z_t],~ m_X=\mathbb{E}\big [X_t],\\
&C_Z = \big (\text{Cov}(Z[k],Z[l])\Big )_{k,l\in [N]},~ C_X = \big (\text{Cov}(X[k],X[l])\Big )_{k,l\in [N]}.
\end{aligned}
\]

From (\ref{mean:000})
\begin{equation}
\begin{aligned}
m_X[k] &= \frac{\varphi - \alpha + \varphib m_Z[k]}
{\alphab + (\varphib/\gamma)(m_Z[k]-\alpha)}.
%&= \begin{cases}
%\frac{\gamma m_Z[k]}{m_Z[k]+ \gamma -\varphi}&\varphi \leq \gamma\\
%\frac{\varphi-\gamma+\varphib m_Z[k]}{\varphi-\gamma+(\varphib/\gamma) m_Z[k]}&\varphi > \gamma
%\end{cases}
\end{aligned}
\label{mean:001}
\end{equation}
The denominator in (\ref{mean:001}) is never zero.
For the covariance matrix,
\begin{equation*}
\begin{aligned}
C_X(k,k) &= m_X[k](1-m_X[k])\\ 
C_X(k,l) &= 
\frac
{\varphi^2A_Z(k,l)(m_X[k]-\gamma)(m_X[l]-\gamma)}
{\alphab^2\gamma^2 - \varphi^2A_Z(k,l)}, \qquad \mbox{for $k \neq l$,}
\end{aligned}
%\label{covarianceconnection:00}
\end{equation*}
where 
\[
A_Z(k,l)=C_Z(k,l) + (m_Z[k]-\alpha)(m_Z[l]-\alpha).
\]
%\begin{Example}[A model only depending on covariances in $C_Z$] 
%For $i,j \in [N]$, $i\ne j$ construct $Z=Z_{ij}$ with probability ${N\choose 2}^{-1}$, where $Z_{ij}[i], Z_{ij}[j]$ have a joint Bernoulli distribution, means $\alpha$, and $Z_{ij}[k]$, $k\ne i,j$ are independent Bernoulli $(\alpha)$ random variables, independent of the $i$ and $j$ th entries. Let $\{Z_t\}$ be time homogeneous with $Z_t =^{\mathcal{D}}Z$.
%The transition distribution is then (\ref{hypercube:000}) where $\kappa_{t,A}$ is zero unless $|A|=2$ and $\kappa_{t,\{i,j\}}= {N\choose 2}^{-1}(\alphab\alpha)^2\text{Cov}\big (Z_{ij}[i],Z_{ij}[j]\big )$.
%%
%\end{Example}
%%
\section{Asymptotic behaviour}\label{Section:limits}
\begin{mdframed}
[style=MyFrame1]
\begin{Proposition}
For $t\in \mathbb{Z}_+$ the coordinates of $X_t$ conditional on $(Z_\tau)_{\tau\leq t}$ are independent and
\begin{equation*}
X_t[k] =^{\cal D} \text{~Bernoulli~}(p_t[k]),
\label{iterate:22a}
\end{equation*}
where 
\begin{equation}
\begin{aligned}
p_t[k]
&= 
\varphi - (\gamma - \varphi)\sum_{l=1}^t
\psi^l\prod_{j=t-l+1}^{t}\Big (1 - Z_j[k]/\alpha\Big )
- \Big ( \varphi - X_{0}[k]\Big )\psi^t
\prod_{j=1}^{t}\Big (1 - Z_j[k]/\alpha\Big ),
\end{aligned}
\label{iterate:22}
\end{equation}

\end{Proposition}
\end{mdframed}
\begin{proof}
The conditional independence follows from Proposition \ref{construction:250}.
For a particular coordinate $k$ denote
\[
\beta := \frac{\varphi}{\gamma},\quad
\rho_{j}[k] :=\frac{\alpha\varphib}{\alphab\varphi}
\Big (1 - \frac{Z_j[k]}{\alpha}\Big ),
\quad j\in \mathbb{Z}_+.
\]
From the versions of (\ref{hypercube:000}) or (\ref{hypercube:001}) conditional on $Z_t$,
\begin{equation}
\begin{aligned}
&\mathbb{E}\Big [1 - X_t[k]/\varphi\ \big | \ X_{t-1}[k],Z_t[k]\Big ]= \rho_{t}[k]\Big (1 - X_{t-1}[k]/\gamma\Big )=\rho_{t}[k]\Bigg (1 - \beta
+\beta\Big (1 - X_{t-1}[k]/\varphi\Big )\Bigg ).
\end{aligned}
\label{iterate:01}
\end{equation}
Iterating in (\ref{iterate:01}) 
\begin{equation}\begin{aligned}
&\mathbb{E}\Big [1 - X_t[k]/\varphi\>\Big | \>X_{0},(Z_\tau)_{\tau\leq t}\Big ]= \sum_{l=1}^t \beta^{l-1}(1-\beta)\prod_{j=t-l+1}^{t}\rho_j[k]
+ \beta^{t}\Big( \prod_{j=1}^{t}\rho_j[k]\Big )\Big ( 1 - X_{0}[k]/\varphi\Big ).
\end{aligned}
\label{iterate:02}
\end{equation}
Noting that 
$
\psi = \frac{\alpha\varphib}{\alphab\varphi}\beta
$
and rearranging (\ref{iterate:02}) gives (\ref{iterate:22}).
The claim is that (\ref{iterate:22}) holds jointly for all $t\in \mathbb{Z}_+$. This follows from the Markov nature of
 $(X_t)_{t}$ conditional on $(Z_t)_{t}$.
\end{proof}
\begin{mdframed}
[style=MyFrame1]
\begin{Corollary}\label{homogeneous:000}
If $(Z_t)_{t}$ is homogeneous in time, then for $t\in \mathbb{Z}_+$ (\ref{iterate:22}) is equivalent in distribution to
\begin{equation}
\begin{aligned}
p_t[k]
&= 
\varphi - (\gamma - \varphi)\sum_{l=1}^t
\psi^l\prod_{j=1}^{l}\Big (1 - Z^\circ_j[k]/\alpha\Big )
- \Big ( \varphi - X_{0}[k]\Big )\psi^t
\prod_{j=1}^{t}\Big (1 - Z^\circ_j[k]/\alpha\Big ),
\end{aligned}
\label{iterate:22b}
\end{equation}
where $(Z^\circ_t)_{t}$ is an \emph{i.i.d.} sequence with $Z^\circ_t =^{\cal D} Z_1$.
Set $\Theta[k]:=\mathbb{E}\big [Z_1[k]\big ]$. We have
\begin{equation*}
\begin{aligned}
\mathbb{P}\big (X_t[k]=1\big ) &= 
\varphi
 + \frac{
 (\gamma-\varphi)\varphib\Theta[k]\big (1 - (- \Theta[k]\varphib/\gamma)^t\big )}
 {\varphi+\gamma-1} - (\varphi - X_0[k]) (-\Theta[k]\varphib/\gamma)^t.
\end{aligned}
\label{marginalmean:0}
\end{equation*}
\end{Corollary}
\end{mdframed}

\begin{mdframed}
[style=MyFrame1]
\begin{Proposition}
Let $(Z_t)_{t}$ be homogeneous in time. Then $(X_t)_{t}$ converges in distribution as $t\to \infty$ to a random variable $Y$, such that the coordinates of $Y$ are independent conditional on an \emph{i.i.d.} sequence $(Z^\circ_t)_{t}$, with $Z^\circ_t =^{\cal D} Z_1$. $Y\in \mathcal{V}_N$ has coordinates conditionally distributed as $N$ Bernoulli random variables with parameters for $k\in [N]$ of 
\begin{equation}
\begin{aligned}
p[k]
&= \varphi - (\gamma - \varphi)\sum_{l=1}^\infty
\psi^l\prod_{j=1}^l(1-Z^\circ_j[k]/\alpha)
=\varphi - (\gamma - \varphi)\sum_{l=1}^\infty
\psi^l
\Big (-\frac{\alphab}{\alpha}\Big )^{S_l^\circ[k]},
\end{aligned}
\label{iterate:22c}
\end{equation}
where $S_l^\circ = \sum_{j=1}^lZ^\circ_j$.
\end{Proposition}
\end{mdframed}
\begin{proof}
Use (\ref{iterate:22b}) in Corollary \ref{homogeneous:000} and let $t\to \infty$. 
Recall that 
$
1 - Z^\circ[k]/\alpha = (-\alphab/\alpha)^{Z^\circ[k]}.
$
We show that $\psi\alphab/\alpha < 1$ so the last term in (\ref{iterate:22b}) tends to zero and the first term converges. Now $\psi\alphab/\alpha = \varphib/\gamma$ in both cases $\varphi \le \gamma$ or $\varphi>\gamma$ and $\varphib/\gamma < 1$ because there is an assumption that $\varphi + \gamma >1$.
\end{proof}
%%
%\begin{Example}[Extreme points]\normalfont
%In the extreme point processes where $Z^\circ_j[k]=z[k]$, $k \in [N]$, $j \in \mathbb{N}$,
%\begin{equation*}
%p[k] = \varphi - \frac{\gamma-\varphi(1 - z[k]/\alpha)}{1 - \psi(1-z[k]/\alpha)}.
%%\label{extremep:00}
%\end{equation*}
%\end{Example}
%%
\begin{Example} [de Finetti Exchangeable coordinates]\label{deFinetti}\normalfont
Assume that $Z_t$ are \emph{i.i.d.} for different $t$ and that for each $t\in \mathbb{Z}_+$ $Z_t[k]$ are exchangeable in $k\in [N]$. Suppose that exchangeability has a de Finetti form so that for any collection $\{\alpha_i\}_{i=1}^r \in \{0,1\}^r$
\begin{equation}
\mathbb{P}\big (Z_t[i_1]=\alpha_1,\cdots,Z_t[i_r] = \alpha_r\big ) = \int_{[0,1]}\theta^{\sum_1^r\alpha_i}\thetab^{r - \sum_1^r\alpha_i}\nu(d\theta).
\label{de:250a}
\end{equation}
Exchangeability is described by (\ref{de:250a}) with $Z_t$ and $X_t$ defined on ${\cal V}_N$, however they can be well defined on ${\cal V}_\infty$ with consistent finite dimensional distributions over $N\in \mathbb{Z}_+$. The index $N$ is supressed, but we think of calculations being for a particular finite $N$-dimensional distribution.
\end{Example}
\begin{mdframed}
[style=MyFrame1]
\begin{Proposition}\label{Proposition:deFinetti}  
Let $Y_\infty$ have the limit distribution of $X_t\in {\cal V}_\infty$ as $t\to \infty$ in a de Finetti model for $Z_t$. The conditional probability of coordinates being equal to one in $X_\infty$  when $\varphi \ne \gamma$ is
\begin{equation}
\begin{aligned}
\mathbb{P}\big (Y_\infty[k]=1 \mid (\theta_j)_{j=1}^\infty\big ) &= \varphi - (\gamma-\varphi)\sum_{l=1}^\infty\psi^l\prod_{j=1}^l(1 - \theta_j/\alpha),
\label{Intro:2}
\end{aligned}
\end{equation}
where $(\theta_j)_{j=1}^\infty$, an \emph{i.i.d.}  sequence of random variables in $[0,1]$.
The unconditional measure is the de Finetti measure of the coordinates of $Y_\infty$. The unconditional probability
\begin{equation*}
\begin{aligned}
\mathbb{P}\big (Y_\infty[k]=1 \big )
&=\varphi - (\gamma - \varphi)\frac{\psi(1-\tilde{\theta}/\alpha)}{1 - \psi(1-\tilde{\theta}/\alpha)}
\end{aligned}
\end{equation*}
where
$
\tilde{\theta} = \int_{[0,1]}\theta\nu(d\theta).
$
\end{Proposition}
\end{mdframed}
If $\varphi=\gamma$ the coordinates of $Y_\infty$ are independent Bernoulli $(\varphi)$.
%\emph{Example plots of $\mathbb{P}\big (Y_\infty[k]=1 \big )$ against $\widetilde{\theta}$ when
%\textcolor{blue}{$\varphi=0.4,\gamma=0.7$}, \textcolor{red}{$\varphi=0.7,\gamma=0.4$}.}
%\begin{center}
%
%\input{curves.tex}
%\end{center}
%

\subsubsection*{Ergodicity}
Let $(X_t)_t$ be a time homogeneous process, driven by \emph{i.i.d} $(Z_j^\circ)_j$. Suppose that $\varphi \ne \gamma$.
Take an approach where the probabilities are random depending on $(Z_j^\circ)_j$.
Let the coordinates of $X_0$  be Bernoulli with parameters $p[k]$, $k\in [N]$,  given by (\ref{iterate:22c}). We want to show that $X_1$ has the same distribution. It is sufficient to show this coordinate wise, because the joint distribution is a mixture over the coordinates of $(Z_j^\circ)_j,Z_1$. Now, recalling that the distribution of $X_0$ is completely determined by  $(Z_j^\circ)_j$,
\begin{equation*}
\begin{aligned}
&\mathbb{P}\big (X_1[k]=1\mid X_0[k],(Z_j^\circ)_j,Z_1)
= \varphi - \varphi\mathbb{E}\big [\big ( 1 - X_1[k]/\varphi\big )\mid X_0[k], (Z_j^\circ)_j,Z_1\big ]\\
&=\varphi - \varphi \cdot \frac{\alpha\varphib}{\alphab\varphi}\big ( 1 - Z_1[k]/\alpha\big )\mathbb{E}\big [\big ( 1 - X_0[k]/\gamma\big )\mid (Z_j^\circ)_j\big ]=\varphi - \gamma\psi\big ( 1 - Z_1[k]/\alpha\big )\big ( 1 - p[k]/\gamma\big )\\
&= \varphi - \gamma\psi\big ( 1 - Z_1[k]/\alpha\big )+ \psi\big ( 1 - Z_1[k]/\alpha\big )\Big (\varphi - (\gamma- \varphi)\sum_{l=1}^\infty
\psi^l\prod_{j=1}^l(1-Z^\circ_j[k]/\alpha)\Big )\\
&=\varphi - (\gamma - \varphi)\cdot\psi\big ( 1 - Z_1[k]/\alpha\big )-(\gamma - \varphi)\cdot\psi\big ( 1 - Z_1[k]/\alpha\big )
\sum_{l=1}^\infty
\psi^l\prod_{j=1}^l(1-Z^\circ_j[k]/\alpha)=^{\cal D}p[k],
\end{aligned}
\end{equation*}
since $Z_1$ is independent of $(Z_j^\circ)_j$ and has the same distribution as $Z$.
The last line follows writing
$Z_1^\#[k] =Z_1[k]$,  $Z_j^\#[k] = Z^\circ_{j-1}[k],\ j \geq 2$, and
\[
p[k]
= \varphi - (\gamma - \varphi)\sum_{l=1}^\infty
\psi^l\prod_{j=1}^l(1-Z^\#_j[k]/\alpha).
\]
Therefore $(X_t)_t$ is ergodic.
%with the limit distribution of $X_t$ as $t \to \infty$ being the same as the stationary distribution.
%%
%%
\section{Exchangeable coordinates}\label{Section:exchangeable}
If $Z_t$ has exchangeable coordinates for $t\in \mathbb{Z}_+$ then $X_t$ has exchangeable
coordinates. The transition distribution of $X_t$ only depends on the Hamming distance $\|x_t\|$, but depends on $x_{t-1}$ through $\|x_{t-1}\|$ and  $\langle x_t,x_{t-1}\rangle$. This is equivalent to the transition distribution only depending on the counts of transitions from zero to either zero or one, and counts of transitions from one to either zero or one. The \emph{p.g.f.} of such transitions is found in Corollary \ref{bigpgf:00}.
Proposition \ref{Hamming:000} shows that the Hamming distance $\|X_t\|$ is Markov and that
the transition density can be expressed as a diagonal expansion of Krawtchouk
polynomials in $\|x_t\|$ and $\|x_{t-1}\|$.
\begin{mdframed}
[style=MyFrame1]
\begin{Proposition}[Properties of the Krawtchouk
polynomials]\label{pr:Kr}
 The Krawtchouk polynomials $\{Q_n(\zeta;N,\alpha)\}_{n=0}^N$ that we use are orthogonal on the Binomial $(N,\alpha)$ distribution, scaled so that $Q_n(0;N,\alpha) = 1$, and
\[
\mathbb{E}\big [Q_n(\zeta;N,\alpha)^2\big ]^{-1} = {N\choose n}\Big (\frac{\alpha}{\alphab}\Big )^n.
\]
Their generating function is
\begin{equation*}
\sum_{n=0}^N{N\choose n}Q_n(\zeta;N,\alpha)s^n = (1-(\alphab/\alpha)s)^\zeta(1+s)^{N-\zeta}.
%\label{Kgen:00}
\end{equation*}
If $\zeta= \sum_{l=1}^N\xi_l$, where $\xi_l\in \{0,1\}$, $l=1,\ldots, N$, then $Q_n(\zeta;N,\alpha)$ is the $n$th elementary symmetric function in $\{1 - \xi_i/\alpha\}_{i=1}^N$. That is
\[
Q_n(\zeta;N,\alpha) = {N\choose n}^{-1}\sum_{\sigma\in S_n}\prod_{l=1}^n(1-\xi_{\sigma(l)}/\alpha),
\]
where $S_n$ is the symmetric group of order $n$.
\end{Proposition}
\end{mdframed}
 For a  proof of Proposition~\ref{pr:Kr}  and a probabilistic description of these polynomials see \cite{DG2012}.
\begin{mdframed}
[style=MyFrame1]
\begin{Proposition}\label{Corollary:exchangeable}
If $Z_t$ has exchangeable coordinates then
\begin{equation}
\kappa_{t,A} := \widetilde{\kappa} _{t,|A|}= \Big (\frac{\alpha}{\alphab}\Big )^{|A|}\mathbb{E}\Big [Q_{|A|}\big ( \|Z_t\|;N,\alpha\big )\Big ].
\label{exchangeable:20}
\end{equation}
The transition distribution is
\begin{equation}
p_{t-1\to t}(\bm{x},\bm{y}) = \varphi^{\|\bm{y}\|}\varphib^{N-\|\bm{y}\|}\Bigg \{
1 + \sum_{k=1}^N{N\choose k}\widetilde{\kappa}_{t,k}
R_k\big (\|\bm{x}\|,\|\bm{y}\|,\langle\bm{x},\bm{y}\rangle\big )
\Bigg \},
\label{hypercube:0001}
\end{equation}
where $R_k$ is the coefficient of ${N\choose k}s^k$ in the generating function
\begin{equation}
G_R(\bm{x},\bm{y};s)=
\big (1+s\big )^{N_{00}}
\big (1 - (\gammab/\gamma)s\big )^{N_{10}}
\big (1 - (\varphib/\varphi)s\big )^{N_{01}}
\big (1 + (\varphib\gammab/\varphi\gamma)s\big )^{N_{11}},
\label{RGF:0}
\end{equation}
with
\[
N_{ab} = |\{k:\bm{x}[k]=a,\bm{y}[k]=b\}|,	\quad a,b\in \{0,1\}.
\]
$\{N_{ab}\}$ has an expression as
\begin{eqnarray*}
N_{00} &=& N - \|\bm{x}\|- \|\bm{y}\| + \langle\bm{x},\bm{y}\rangle \nonumber,  
 \qquad  N_{01} = \|\bm{y}\| - \langle\bm{x},\bm{y}\rangle\nonumber, \\
N_{10} &=& \|\bm{x}\| - \langle\bm{x},\bm{y}\rangle\nonumber,  \qquad\qquad \qquad \qquad
N_{11} = \langle\bm{x},\bm{y}\rangle.
\label{Nxy:00}
\end{eqnarray*}
%Transitions $\bm{x}\to \bm{y}$ where coordinates are exchangeably distributed with configuration $\{N_{ij}\}_{i,j=0}^1$ have equal probability.
\end{Proposition}
\end{mdframed}
\begin{proof}
If $Z_t$ has exchangeable coordinates then 
\begin{equation*}
p_{t-1\to t}(\bm{x},\bm{y}) = \varphi^{\|\bm{y}\|}\varphib^{N-\|\bm{y}\|}\Bigg \{
1 + \sum_{k=1}^N{N\choose k}\widetilde{\kappa}_{t,k}
R_k^\circ(\bm{x},\bm{y})
\Bigg \},
\end{equation*}
where 
\[
R_k^\circ(\bm{x},\bm{y})
= {N\choose k}^{-1}\sum_{A\subseteq[N], |A|=k}
\prod_{j\in A}\Big (1 - \frac{x[j]}{\gamma}\Big )
\Big (1 - \frac{y[j]}{\varphi}\Big ).
\]
Now
\[
\Big (1 - \frac{x[j]}{\gamma}\Big )
\Big (1 - \frac{y[j]}{\varphi}\Big )=
\begin{cases}
1&\text{if~}x[j]=0,y[j]=0\\
-\gammab/\gamma&\text{if~}x[j]=1,y[j]=0\\
-\varphib/\varphi&\text{if~}x[j]=0,y[j]=1\\
\gammab\varphib/\gamma\varphi&\text{if~}x[j]=1,y[j]=1
\end{cases},
\]
so considering the counts $\{N_{ab}\}$, $a,b\in \{0,1\}$, we have 
$
R_k^\circ(\bm{x},\bm{y})=R_k\big (\|\bm{x}\|,\|\bm{y}\|,\langle\bm{x},\bm{y} \rangle\big )
$, 
with generating function (\ref{RGF:0}).
\end{proof}
If $\varphi=\gamma$ then Corollary \ref{Corollary:exchangeable} is equivalent to Lemma 4 in \cite{CG2021}, however the expression in terms of $R_k$ is new.

\begin{mdframed}
[style=MyFrame1]
\begin{Proposition}\label{Hamming:000}
$\|X_t\|$ is the Hamming distance of $X_t$ from $\bm{0}$.
If $Z_t$ has exchangeable coordinates then $(\|X_t\|)_{t}$ is a Markov chain with transition distribution 
\begin{equation}
\begin{aligned}
&p_{t-1\to t}(\|\bm{y}\|\mid \|\bm{x}\|)
= {N\choose \|\bm{y}\|}\varphi^{\|\bm{y}\|}\varphib^{N-\|\bm{y}\|}\Big \{1 + \sum_{k=1}^N{N\choose k}\widetilde{\kappa}_{t,k}
Q_k(\|\bm{y}\|;N,\varphi)Q_k( \|\bm{x}\|;N,\gamma)\Big \}.
\end{aligned}
\label{Hammingtransition:00}
\end{equation}
\end{Proposition}
\end{mdframed}
\begin{proof}
To show that $(\|X_t\|)_{t}$ is Markov note that for $A\subseteq [N]$ with 
$|A|=k$
\[
\mathbb{E}\Big [\prod_{j\in A} \Big ( 1 - \frac{\bm{Y}[j]}{\varphi}\Big )\mid \|\bm{y}\|\Big ]
= Q_k(\|\bm{y}\|;N,\varphi).
\]
This implies (\ref{Hammingtransition:00}) because 
\begin{equation*}
\begin{aligned}
&p_{t-1\to t}(\|\bm{y}\|\mid \bm{x})={\rm Bin}(N, \varphi, \|\bm{y}\|)
\Big \{1 + \sum_{k=1}^N \widetilde{\kappa}_{t,k}
Q_k(\|\bm{y}\|;N,\varphi)
\sum_{A\subseteq [N], |A| = k}\prod_{j\in A} \Big ( 1 - \frac{\bm{x}[j]}{\varphi}\Big )\Big \}\\
&~={\rm Bin}(N, \varphi, \|\bm{y}\|)
\Big \{1 + \sum_{k=1}^N \widetilde{\kappa}_{t,k}
{N\choose k}Q_k(\|\bm{y}\|;N,\varphi)Q_k(\|\bm{x}\|;N,\gamma)\Big \}.
\end{aligned}
\end{equation*}
\end{proof}
If $Z_t$ has does not have exchangeable coordinates then $(\|X\|_t)_{t}$ is not  Markovian because $\kappa_{t,A}$ can depend on $A$, not only on $|A|$.
 \cite{DG2012} have a characterization of conditional bivariate Binomial distributions which have diagonal expansions in Krawtchouk polynomials with two different parameters. The class of transition distributions $p_{t-1\to t}(\|\bm{y}\|,\|\bm{x}\|)$ with $\widetilde{\kappa}_{t,k}$ taking the form
 (\ref{exchangeable:20}) is identical to their characterized class.% The extreme point chains homogeneous in time have
\begin{Example} [Extreme point chains]\label{Example:ExtremePointChains}\normalfont
There is no simple general solution to finding
% the sequence  $\{c_r\}$ in
  the distribution
 % (\ref{stationary:101}) 
  of the Hamming distance, however we can make progress considering the extreme points as $N\to \infty$.
Let $\|Z_t\|=M\in [N]$, a fixed number for $t\in \mathbb{N}$, such that the $M$ unit entries in $Z_t$ are uniformly distributed for each $t\in \mathbb{N}$. Then
$
\kappa_r = (\alpha/\alphab)^rQ_r(M;N,\alpha),\> r \in [N]
$.

As $N\to \infty$ then the processes at each coordinate are effectively independent in the limit taking $Z_t$ to have a product Bernoulli$(M/N)$ distribution. The stationary distribution of $\|X_t\|$ is asymptotically Binomial$(N,p)$ with 
\begin{equation*}
\begin{aligned}
p &= \frac
{
\varphi - \frac{\varphib\alpha}{\alphab}
\Big ( 1 -\frac{1}{\alpha}(M/N)\Big )
}
{
1 - \frac{\alpha\varphib}{\alphab\gamma}\Big ( 1 -\frac{1}{\alpha}(M/N)\Big )
}
=\begin{cases}
\frac{M/N}
{1 - \varphi/\gamma}+\mathcal{O}(N^{-2}),
&\varphi<\gamma\\
1 - \frac{(\varphib/\gamma)(M/N)}
{1 - \varphib/\gammab} +\mathcal{O}(N^{-2}),
&\varphi>\gamma.
\end{cases}
\end{aligned}
\label{Bernp:00a}
\end{equation*}
from (\ref{Bernp:00}). In the two cases for $\alpha$, suppose $M$ is taken as fixed. In the limit as $N\to \infty$, and $\varphi<\gamma$, then the stationary distribution of $\|X_t\|$ is Poisson $\big (M/(1-\varphi/\gamma)\big )$. If $\varphi > \gamma$ then the stationary distribution of $N- \|X_t\|$ is Poisson $\big ((\varphib/\gamma) M/(1 - \varphib/\gammab)\big )$.

There is another limit where $M\to \infty$ such that $M/N\to \omega\in (0,1)$. Then 
\[
p \to \frac
{
\varphi - \frac{\varphib\alpha}{\alphab}
\Big ( 1 -\frac{1}{\alpha}\omega\Big )
}
{
1 - \frac{\alpha\varphib}{\alphab\gamma}\Big ( 1 -\frac{1}{\alpha}\omega\Big )
}
=\begin{cases}
\frac{\omega}
{1 - \varphi/\gamma},%+\mathcal{O}(N^{-2}),
&\varphi<\gamma\\
1 - \frac{(\varphib/\gamma)\omega}
{1 - \varphib/\gammab}
\end{cases}
\]
which is not zero or one. Then the stationary distribution of the standardized value of $\|X_t\|$ converges to $N(0,1)$.
\end{Example}
\begin{mdframed}
[style=MyFrame1]
\begin{Corollary}\label{bigpgf:00}
Let $(X_t)_t$ be  time-homogeneous with exchangeable coordinates.
The pgf of $(N_{01},N_{11})$, conditional on  $\|Z_t\|=\zeta$ is the coefficient of  ${N\choose \zeta}\xi^\zeta$, where $N = N_{01} +N_{11}$,
in 
\begin{equation}
\begin{aligned}
&
\begin{cases}
(1+\xi s_{01})^{N_{0\cdot}}
\Big (1 - \varphi/\gamma + s_{11}\varphi/\gamma
+ \xi(\varphib/\gamma + s_{11}(1-\varphib)/\gamma )\Big )^{N_{1\cdot}},
&\varphi \leq \gamma\\
\big (\varphib/\gammab + s_{01}(1-\varphib/\gammab ) + \xi s_{01}\big )^{N_{0\cdot}}
\big ( s_{11} + \xi(\varphib/\gamma + s_{11}(1-\varphib/\gamma)\big )^{N_{1\cdot}},
&\varphi > \gamma
\end{cases}
.
\end{aligned}
\label{pgf2:00}
\end{equation}
\end{Corollary}
\end{mdframed}
\begin{proof}
 Note that
\[
\widetilde{\kappa}_{t,k} 
= \Big (\frac{\alpha}{\alphab}\Big )^kQ_k(\zeta;N,\alpha)
= \Big (\frac{\alpha}{\alphab}\Big )^kQ_\zeta(k;N,\alpha),
\]
which is the coefficient of ${N\choose \zeta}\xi^\zeta$ in 
$
(\alpha/\alphab)^k(1 - \xi\alphab/\alpha)^k( 1 + \xi)^{N-k}
$.\\
From (\ref{hypercube:0001}) and (\ref{RGF:0}) $p_{t-1\to t}(\bm{x},\bm{y})$ is the coefficient of ${N\choose \zeta}\xi^\zeta$ in
$
(1 + (\alpha/\alphab)\xi)^NG_R(\bm{x},\bm{y};s^*)
$,
where $s^*=(\alpha/\alphab)(1-\xi\alphab/\alpha)/(1 + \xi)$.
The joint \emph{p.g.f.} of $N_{01},N_{11}$ conditional on $N_{0\cdot},N_{1,\cdot}$ is the coefficient of ${N\choose \zeta}\xi^\zeta$ in $( 1 + \xi)^N$ times
\[
\mathbb{E}\big [s_{01}^{W_{01}}s_{11}^{W_{11}}
\big (1+s^*\big )^{N_{0\cdot} - W_{01}}
\big (1 - (\varphib/\varphi)s^*\big )^{W_{01}}
\big (1 + (\varphib\gammab/\varphi\gamma)s^*\big )^{W_{11}}
\big (1 - (\gammab/\gamma)s^*\big )^{N_{1\cdot}-W_{11}} \mid N_{0\cdot},N_{1\cdot}\big ]
\]
where $W_{01}, W_{11}$ are independent  Binomial with parameters $(N_{0\cdot},\varphi),(N_{1\cdot},\varphi)$.
That is
\[
\begin{aligned}
&\mathbb{E}\big [s_{01}^{N_{01}}s_{11}^{N_{11}}\mid N_{0\cdot},N_{1\cdot}\big ]\\
&=(1 + \xi)^N\big (\varphib(1+s^*) + s_{01}\varphi(1-(\varphib/\varphi)s^*)\big )^{N_{0\cdot}}
\big (\varphib ( 1 - (\gammab/\gamma)s^*) + s_{11}\varphi(1+ (\varphib\gammab)/(\varphi\gamma)s^*)\big )^{N_{1\cdot}}\\
%&=(1 + \xi)^N\big (\varphib+\varphi s_{01} + s^*\varphib(1-s_{01})\big )^{N_{0\cdot}}
%\big (\varphib+\varphi s_{11} - s^*(\varphib\gammab/\gamma)(1-s_{11})\big )^{N_{1\cdot}}\\
%&=\big ((1 + \xi)(\varphib+\varphi s_{01}) + (\alpha/\alphab)(1-\xi\alphab/\alpha)\varphib(1-s_{01})\big )^{N_{0\cdot}}\\
%&~~\times\big ((1 + \xi)(\varphib+\varphi s_{11}) - (\alpha/\alphab)(1-\xi\alphab/\alpha)(\varphib\gammab/\gamma)(1-s_{11})\big )^{N_{1\cdot}}\\
&= \big (\varphib/\alphab + s_{01}(1-\varphib/\alphab) + s_{01}\xi\big )^{N_{0\cdot}}\\
&~~\times
\big (\varphib(1-(\alpha/\alphab)(\gammab/\gamma)) + s_{11}(1-(\varphib(1-(\alpha/\alphab)(\gammab/\gamma) )
+ \xi (\varphib/\gamma + s_{11}(1-\varphib/\gamma)\Big )^{N_{1\cdot}}.
\end{aligned}
\]
Simplifying gives (\ref{pgf2:00}).
\end{proof}
\begin{mdframed}
[style=MyFrame1]
\begin{Corollary}\label{jointpgf:600}
Assuming the conditions in Corollary \ref{bigpgf:00} hold and in addition $\zeta$ has a Binomial $(N,\theta)$ distribution then the joint \emph{p.g.f.}
\begin{equation}
\begin{aligned}
&\mathbb{E}\big [s_{01}^{N_{01}}s_{11}^{N_{11}}h^\zeta\mid N_{0\cdot},N_{1\cdot}\big ]\\
&
=\begin{cases}
(\thetab+\theta h s_{01})^{N_{0\cdot}}
\Big (\thetab(1 - \varphi/\gamma + s_{11}\varphi/\gamma)
+ \theta h(\varphib/\gamma + s_{11}(1-\varphib/\gamma) )\Big )^{N_{1\cdot}},
&\varphi \leq \gamma\\
\big (\thetab(\varphib/\gammab + s_{01}(1-\varphib/\gammab )) + \theta h s_{01}\big )^{N_{0\cdot}}
\big ( \thetab s_{11} + \theta h(\varphib/\gamma + s_{11}(1-\varphib/\gamma)\big )^{N_{1\cdot}},
&\varphi > \gamma
\end{cases}.
\end{aligned}
\label{pgf2:01}
\end{equation}
$N_{01},N_{11}$ are independent Binomial with respective parameters
\begin{equation}
\begin{aligned}
&(N_{0\cdot}, 1 - \thetab\varphib/\alphab),\quad
(N_{1\cdot}, (\alpha\thetab + \theta(\varphi+\gamma-1))/\gamma)\\
&\equiv
\begin{cases} 
(N_{0\cdot}, \theta),\quad
(N_{1\cdot}, (\varphi - \theta\gammab)/\gamma),& \varphi \leq \gamma\\
(N_{0\cdot},1 - \thetab\varphib/\gammab)\quad (N_{1\cdot},1 - \theta\varphib/\gamma),&\varphi > \gamma
\end{cases}.
\label{binomials:02}
\end{aligned}
\end{equation}
The covariances between $\zeta$ and $N_{01},N_{11}$ are
\begin{equation*}
\begin{aligned}
&\text{Cov}(\zeta,N_{01})  = N_{0\cdot}\theta\thetab\varphib/\alphab,
\quad \text{Cov}(\zeta,N_{11}) = -N_{1\cdot}\theta\thetab\varphib\gammab/(\gamma\alphab)\\
&\equiv\begin{cases}
\text{Cov}(\zeta,N_{01})  = N_{0\cdot}\theta\thetab,
\quad \text{Cov}(\zeta,N_{11}) = -N_{1\cdot}\theta\thetab\gammab/\gamma,&\varphi \leq \gamma\\
\text{Cov}(\zeta,N_{01})  = N_{0\cdot}\theta\thetab\varphib/\gammab,
\quad \text{Cov}(\zeta,N_{11}) = -N_{1\cdot}\theta\thetab\varphib/\gamma,&\varphi > \gamma
\end{cases}.
\end{aligned}
\label{Bcovars:0}
\end{equation*}
\end{Corollary}
\end{mdframed}
\begin{proof}
All that is needed is to replace $\xi$ by $h\theta/\thetab$ in (\ref{pgf2:00}) and multiply by $\thetab^N$.
Setting $h=1$ in (\ref{pgf2:01}) shows the Binomial structure. The probabilities in the Binomial parameters  \eqref{binomials:02} must be in $[0,1]$ and a check shows that this is true. The means and covariances are straightforward to calculate from (\ref{pgf2:01}).
\end{proof}
\section{Normal limit theorems as $N\to \infty$}\label{Section:normal}
Limit distributions as $N\to \infty$ for the Hamming distance are considerably simpler than the exact discrete distributions.

In an exchangeable model the transition counts $\{N_{ij}\}_{i,j=0}^1$ have a Central limit theorem that arises from the Binomial counts in Corollary \ref{bigpgf:00} and \ref{jointpgf:600}. 
\begin{mdframed}
[style=MyFrame1]
\begin{Proposition}\label{CentralN:00}
Let $(Z_t^{(N)})_t$ be homogeneous in time with exchangeable coordinates, and suppose that $\lim_N \|Z_t^{(N)}\|/N \to \Theta_t$, a random variable in $[0,1]$.   We have that
\begin{equation*}
V_{01} = \frac{N_{01} - N_{0\cdot}q_{01}}{N_{0\cdot}^{1/2}},\quad
V_{11} = \frac{N_{11} - N_{1\cdot}q_{11}}{N_{1\cdot}^{1/2}},
\end{equation*}
conditional on $N_{0\cdot},N_{1\cdot},\Theta_t=\theta$ as $N_{0\cdot},N_{1\cdot} \to \infty$, converge to independent Normal with
with zero means, and variances $q_{01}\bar{q}_{01}$, $q_{11}\bar{q}_{11}$,
 where
$q_{01}=1 - \thetab\varphib/\alphab$, $q_{11} = (\alpha\thetab + \theta(\varphi+\gamma-1))/\gamma$.
\end{Proposition}
\end{mdframed}
\begin{proof}
The Binomial type form of the \emph{p.g.f.} (\ref{pgf2:01}) implies the Central limit theorem. We omit the detailed calculation. 
It is enough to use the means, variances and covariances calculated in Corollary \ref{jointpgf:600}. 
In the  \emph{p.g.f.} an assumption was that $\|Z_t^{(N)}\|$ is Binomial $(N,\theta)$. Calculations for this Proposition  assume that 
$\|Z_t^{(N)}\|/N  =  \Theta_t$ and the arguement is then condtional on $\Theta_t=\theta$. This is the same in the limit as taking $\|Z_t^{(N)}\|$ to be Binomial $(N,\theta)$ because then $\|Z_t^{(N)}\|/N \to \theta$.
\end{proof}
\begin{mdframed}
[style=MyFrame1]
\begin{Proposition}\label{HammingNormal:00}
Let $(Z_t^{(N)})$ be homogeneous in time with exchangeable coordinates, and suppose that $\lim_N \|Z_t^{(N)}\|/N = \Theta_t$, a random variable in $[0,1]$.  Let $V_t,\ U_t$ be distributed as the limits of 
\begin{equation*}
\begin{aligned}
V_t^{(N)} = N^{-1/2}(\|X_t\| - N\varphi), \qquad \mbox{and} \qquad 
U_t^{(N)} =  N^{-1/2}(\|X_{t-1}\| - N\gamma)
\end{aligned}
\end{equation*}
as $N\to \infty$. $V_t$ has a Normal distribution, conditional on $\Theta_t=\theta,\ U_t=u$, 
with mean and variance
\begin{equation*}
\begin{aligned}
\mu_\infty = \psi\big (1-\theta/\alpha\big )u \qquad \mbox{and} \qquad 
\sigma^2_\infty = 
\varphi\varphib\Big (1 - (\varphib/\varphi)(\gammab/\gamma)(\alpha/\alphab)^2(1-\theta/\alpha)^2\Big ).
\end{aligned}
%\label{mv:256}
\end{equation*}
%If $\Theta_t=\theta$, a constant for all $t$ then
%the stationary distribution is normal with zero mean and variance
%\begin{equation}
%\frac{\sigma^2_\infty}{1 - \psi^2(1-\theta/\alpha)^2}.
%\label{var:620}
%\end{equation}
\end{Proposition}
\end{mdframed}
\begin{proof}
%Let $N_{ij}=|\{k:X^{(N)}_{t-1}[k]=i,\ X_t^{(N)}[k]=j\}|$, for $i,j=0,1$.
Note that $\|X_t^{(N)}\|=N_{01}+N_{11}$ and $\|X_{t-1}^{(N)}\|=N_{1\cdot}$.
\begin{equation}
\begin{aligned}
&N^{-1/2}\Big (\|X_t^{(N)}\| - N\varphi\Big ) =N^{-1/2}\Big ( N_{01} + N_{11} - N\varphi\Big )\\
%& =N^{-1/2}\Big ( N_{01} - N_{0\cdot}q_{01} + N_{11} - N_{1\cdot}q_{11}
% + N_{0\cdot}q_{01} + N_{1\cdot}q_{11}  -N\varphi\Big )\\
%&= N^{-1/2}\Big (N_{01} - N_{0\cdot}q_{01} + N_{11} - N_{1\cdot}q_{11}
% + (N_{1\cdot}-N\gamma)(q_{11}-q_{01})\Big )+ N^{1/2}\big (q_{01} + \gamma(q_{11}-q_{01}) - \varphi\big )\\
&=
(N_{0\cdot}/N)^{1/2}V_{01} + (N_{1\cdot}/N)^{1/2}V_{11} + N^{-1/2}  (N_{1\cdot}-N\gamma)(q_{11}-q_{01})+ N^{1/2}\big (q_{01} + \gamma(q_{11}-q_{01}) - \varphi\big ).
\end{aligned}
\label{decompose:25}
\end{equation}
In last term of (\ref{decompose:25}), $q_{01} + \gamma(q_{11}-q_{01}) - \varphi=0$ and $N_{0\cdot}/N \to \gammab$,
$N_{1\cdot}/N \to \gamma$. Therefore $N^{-1/2}\Big (\|X_t^{(N)}\| - N\varphi\Big )$ converges in distribution to
a normal random variable with mean and variance
\[
\begin{aligned}
&\mu_\infty = (q_{11}-q_{01})u
= \psi\big (1-\theta/\alpha\big )u,\\
&\sigma^2_\infty = \gammab q_{01}\bar{q}_{01}+\gamma q_{11}\bar{q}_{11}
= \varphi\varphib\Big (1 - (\varphib/\varphi)(\gammab/\gamma)(\alpha/\alphab)^2(1-\theta/\alpha)^2\Big ).
\end{aligned}
\]
%We calculate the stationary variance when $\Theta_t=\theta$ for all $t$.
%\[
%\begin{aligned}
%\text{Var}\big (V) &= \mathbb{E}\big [\text{Var}\big (V\mid U\big )\big ]
%+ \text{Var}\big (\mathbb{E}\big [V\mid U\big ]\big )= \sigma^2_\infty + \psi^2\big(1 - \theta/\alpha)^2\text{Var}\big (U\big ).
%\end{aligned}
%\]
%Equating $\text{Var}\big (V\big )=\text{Var}\big (U\big )$,
%$\text{Var}\big (V\big )$ is seen to be equal to (\ref{var:620}).
% The stationary mean is clearly zero.
\end{proof}
%%
%A direct proof of convergence of the transition density (\ref{HammingNormal:01}) to a Normal density is possible using a property that the Krawtchouk polynomials converge to Hermite polynomials. Details are not included in this paper.
%%%
% %%
Proposition \ref{HammingNormal:00} finds the normal limit of the transition distribution and the stationary distribution in the limit transition distribution.
Finding the limit distribution of the process $N^{-1/2} \big (\|X_t^{(N)}\| -N\varphi\big )$ in the next Proposition \ref{HammingNormal:01} is different because of the identical scaling $N^{-1/2} \big (\|X_{t-1}^{(N)}\| -N\varphi\big )$,
subtracting the mean $N\varphi$. To obtain such a limit what is required for stability is that $\varphi,\gamma$ depend on $N$ and $N^{1/2}(\varphi^{(N)}-\gamma^{(N)}) \to c$, a constant. It is convenient to take $\alpha = \min \{\varphi^{(N)},\gamma^{(N)}\}$ not depending on $N$. 
% The limit process is expressed as an AR(1) series. 
%% New version
\begin{mdframed}
[style=MyFrame1]
\begin{Proposition}\label{HammingNormal:01}
Let $(Z_t^{(N)})$ be homogeneous in time with exchangeable coordinates, and suppose that $\lim_N \|Z_t^{(N)}\|/N = \Theta_t$, a random variable in $[0,1]$ with probability measure $\nu$.  Assume that $\varphi^{(N)},\gamma^{(N)}$ depend on $N$ with
 $N^{1/2}(\varphi^{(N)}-\gamma^{(N)}) \to c$, a constant. $\alpha$ is kept fixed as $N\to \infty$. Let $V_t$ be distributed as the limit of 
$
V_t^{(N)} = N^{-1/2}(\|X_t\| - N\varphi^{(N)})
$
as $N\to \infty$. 
$V_t$ is normal  
with mean and variance, conditional on $\Theta_t$,
\begin{equation}
\mu_t = \big (1-\Theta_t/\varphi\big )(v_{t-1}+c),\quad
\sigma^2_t = 
\varphi\varphib\Big (1 - (1-\Theta_t/\varphi)^2\Big ).
\label{mv:250}
\end{equation}
An AR(1) series representation with random parameters is
$
V_t = (1-\Theta_t/\varphi)\big (c + V_{t-1}\big ) + \varepsilon_t,
$
where 
$\varepsilon_t$ is independent of $V_{t-1}$ with 
$\text{Var}(\varepsilon_t) = \varphi\varphib\big (1 - (1-\Theta_t/\varphi)^2\big )$.

If $\Theta_t=\theta$, a constant, then $V_t$ is a normal AR(1) series with stationary mean $\varphi c(1-\theta/\varphi)/\theta$ and variance $\varphi\varphib$.
\end{Proposition} 
\end{mdframed}
\begin{proof}
The proof is a modification of the proof in Proposition \ref{HammingNormal:00}, conditioning on $\Theta_t=\theta$, when 
\[
\begin{aligned}
v_{t-1}&=N^{-1/2} \big (\|\bm{x}^{(N)}\| - N\varphi^{(N)}\big )= N^{-1/2} \big (\|\bm{x}^{(N)}\| - N\gamma^{(N)}\big )-N^{1/2}(\varphi^{(N)} - \gamma^{(N)})\to u - c.
\end{aligned}
\]
$\mu_t$ and $\sigma_t$ are modifications of $\mu_\infty$ and $\sigma^2_\infty$ obtained by setting $\varphi=\gamma=\alpha$, in their limit values and taking into account that $u=v_{t-1}+c$.

To calculate the stationary distribution of $V_t$, when $\Theta_t = \theta$, a constant for all $t$, let $V_{t-1}$ have a N$(\delta,\beta^2)$ distribution. Then $V_t$ has a Normal distribution with mean $(1-\theta/\varphi)(\delta + c)$
and variance $\sigma^2_t + (1-\omega/\varphi)^2\beta^2$. Solving $\delta = (1-\theta/\varphi)( c+\nu)$ and
$\beta^2 = \sigma^2_t + (1-\theta/\varphi)^2\beta^2$, the stationary mean and variance are
$\delta = \varphi c(1-\theta/\alpha)/\theta$ and $\beta^2=\sigma_t^2/\big (1 - (1-\theta/\varphi)^2\big )=\varphi\varphib$.
\end{proof} 
%%%
The Hermite polynomials $\big (H_k(\cdot;\sigma^2)\big)_{n=0}^\infty$ are used in the next Corollary. They are  orthogonal on the N$(0,\sigma^2)$ distribution and have a generating function
\begin{equation*}
\sum_{k=0}^\infty \frac{z^k}{k!}H_k(w;\sigma^2) = \exp \big \{zw - \frac{1}{2}\sigma^2z^2\Big \}.
%\label{HermiteGF:00}
\end{equation*}
Note that $H_k(w;\sigma^2) = \sigma^kH_k(w/\sigma;1)$. 
Their transforms are 
\[
\mathbb{E}\big [e^{\vartheta W}H_j(W;\sigma^2)\big ]= e^{\sigma^2\vartheta^2/2}\big (\vartheta\sigma^2\big )^j/j!.
\]
\begin{mdframed}
[style=MyFrame1]
\begin{Corollary}
Let $\Phi_t= \big ( 1 - \Theta_t/\varphi\big )$. The transition density of $V_t$ given $V_{t-1}=v_{t-1}$ is
\begin{equation}
\frac{1}{\sqrt{2\pi\varphi\varphib}}\exp\big \{\frac{-v_t^2}{\varphi\varphib}\big  \}
\Big \{1 + \sum_{j=1}^\infty\frac{\mathbb{E}\big [\Phi_t^j\big ]}{j!}H_j(v_t/\sqrt{\varphi\varphib};1)H_j((v_{t-1}+c)/\sqrt{\varphi\varphib};1)\Big \}.
\label{density:900}
\end{equation}
\end{Corollary}
\end{mdframed}
\begin{proof}

The \emph{m.g.f.} of $V_t$ given $V_{t-1} = v_{t-1}$, using (\ref{mv:250}) is
\begin{equation}
e^{\vartheta\Phi_t(c+v_{t-1}) +\vartheta^2 \varphi\varphib\big ( 1 - \Phi_t^2\big )/2}
= e^{-\vartheta^2 \varphi\varphib}\sum_{j=0}^\infty \frac{\vartheta^j\Phi_t^j}{j!}H_j(v_{t-1}+c;\varphi\varphib).
\label{mgf:900}
\end{equation}
Inverting (\ref{mgf:900}) gives (\ref{density:900}).

\end{proof}

A time dependent solution for the mean in the AR(1) model when $\Theta_t=\theta$ for all $t$ is
\[
\mathbb{E}\big [V_t\big ] 
= \varphi c(1-\theta/\varphi)/\theta + (1-\theta/\varphi)^t\Big (\mathbb{E}\big [V_0\big] - \varphi c(1-\theta/\varphi)/\theta\Big )
\]
with the rate of approach $(1-\theta/\varphi
)^t$ as $t\to \infty$.  Since $\theta \in (0,1)$,
$1-\theta/\varphi \in (-\varphib/\varphi, 1)$. The absolute value of the lower terminal is less than or equal to 1 because $\varphi \geq 1/2$ under our assumption that $\varphi \geq 1/2$ as a limit from $\varphi^{(N)}+\gamma^{(N)} \geq 1$.

Proposition \ref{HammingNormal:01} shows that we must take the scaling
$N^{1/2}(\varphi^{(N)}-\gamma^{(N)}) \to c \in (-\infty,\infty)$ as $N\to \infty$ to obtain a proper AR(1) series. If $\varphi, \gamma$ do not depend on $N$ and $\varphi\ne\gamma$, then $\mathbb{E}\big [V_t\big] \to \pm\infty$, with the sign equal to $\text{Sgn} \big ((1-\omega/\alpha)(\varphi-\gamma)\big )$.
\bigskip

In Propositions \ref{CentralN:00}, \ref{HammingNormal:00} and  \ref{HammingNormal:01} our interest is in Central limit theorems for the Hamming distance, however another statement is 
\begin{equation}
\frac{N_{01}}{N_{0\cdot}} \to q_{01},\ \frac{N_{11}}{N_{1\cdot} }\to q_{11}
\label{weak:00}
\end{equation}
as $N_{0\cdot},\ N_{1\cdot} \to \infty$ conditional on $N_{0\cdot}, N_{1\cdot}, \Theta_t=\theta$. A process $W_t$ approximating $\|X_t\|$ for large $N$, based on (\ref{weak:00}) is therefore
\begin{equation}
\begin{aligned}
W_t &= \big (N - W_{t-1}\big )\big ( 1 - \bar{\Theta}_t\varphib/\alphab\big )
+ W_{t-1}\big (\alpha \bar{\Theta}_t + \Theta_t(\varphi+\gamma - 1\big )/\alpha\\
&= N\big ( 1 - \bar{\Theta}_t\varphib/\alphab\big ) + W_{t-1}\psi (1-\Theta_t/\alpha\big ).
\end{aligned}
\label{ApproxN:00}
\end{equation}
It is easy to simulate a sample path $(W_t)$ from (\ref{ApproxN:00}).

 \section{Estimation from a sample path}
 Let $(X_t)$ be time homogeneous and assume that the coordinates of $X_t$ are exchangeable.
The data is a sample path $(x_t)_{t=0}^T$. Let $n_{ij}^t$, $i,j \in \{0,1\}$ denote the number of transitions $i \to j$ in transition $t$. 
\subsubsection*{Estimation of $\varphi,\gamma$}
If the coordinates of $X_{t-1}$ have an independent $N$-Bernoulli$(\gamma)$ distribution then the coordinates of $X_t$ have an independent $N$-Bernoulli$(\varphi)$ distribution. This is always true and does not depend on $(Z_t)$.
The two parameters can be estimated by a least squares approach by minimizing
\begin{equation}
\sum_{t=1}^T \Big (\gammab \widehat{p}_{01}[t] + \gamma \widehat{p}_{11}[t] - \varphi\Big )^2,
\label{ls:10}
\end{equation}
where $\widehat{p}_{ij}[t] = n_{ij}^t/n^{t-1}_{i\cdot}$, assuming that $n^{t-1}_{i\cdot}\geq 1$. Denote
$d[t] = \widehat{p}_{01}[t]-\widehat{p}_{11}[t]$, $\bar{p}_{ij} = \sum_{t=1}^T\widehat{p}_{ij}[t]/T$, $\bar{d} = \sum_{t=1}^Td[t]/T$.
Suppose $\varphi \ne \gamma$. The least squares estimates satisfy
\begin{equation}
\begin{aligned}
\widehat{\varphi} + \widehat{\gamma}\bar{d}= \bar{p}_{01}, \qquad \mbox{ and } \qquad 
\widehat{\varphi}\bar{d}
+ \widehat{\gamma}\sum_{t=1}^Td[t]^2/T = 
\sum_{t=1}^Td[t]\widehat{p}_{01}[t]/T.
\end{aligned}
\label{regression:250}
\end{equation}
Solving (\ref{regression:250})
\begin{equation*}
\widehat{\gamma} = 
\frac{
\sum_{t=1}^T(d[t]-\bar{d})(\widehat{p}_{01}[t]-\bar{p}_{01})
}
{
\sum_{t=1}^T(d[t]-\bar{d})^2
},
\end{equation*}
with $\widehat{\varphi}$ determined by the first equation in (\ref{regression:250}).\\
If $\varphi = \gamma$ minimizing (\ref{ls:10}),
$
\widehat{\varphi} = \frac{\bar{p}_{01}}{1+\bar{d}}.
$

 \subsection*{Estimating the distribution of $Z_t$}\label{Section:estimation}
Estimation of the distribution of $Z_t$ in a simple exchangeable model is now worked through.  We assume that $\varphi$ and $\gamma$ are known, but if not they can be estimated by the proceedure above. At each transition, conditional on an observation $\theta_t$ from a measure $\nu$ take the $N$ coordinates of $Z_t$ to be \emph{i.i.d.} Bernoulli $(\theta_t)$.  We assume that $\nu$ is continuous with support a subset or equal to $(0,1)$ and a density with respect to Lebesgue measure. $X_t$ is observed but $Z_t$ is not.
This is a Bayesian approach where in each transition the probability in the Bernoulli distribution of coordinates is a random variable.
  The likelihood in transition $t$, conditional on $\theta_t$, when $\{n_{ij}^t\}$ is observed is
\begin{equation*}
\prod_{i=0}^1\prod_{j=0}^1{n_{i\cdot}\choose n_{ij}}p_t(i,j;\theta_t)^{n_{ij}^t},
%\label{like:342}
\end{equation*}
where
\begin{equation}
p_t(i,j;\theta_t) = 
\varphi^j\varphib^{1-j}
\Big \{1 + \frac{\alpha}{\alphab}\big (1 - \frac{\theta_t}{\alpha} \big )
 \big (-\frac{\varphib}{\varphi}\big )^{j}
\big (-\frac{\gammab}{\gamma}\big )^{i}\Big )
\Big \}
\label{like:343}
\end{equation}
because each transition from $i\to j$ has the same probability in different coordinates.
The probabilities in (\ref{like:343}) are from (\ref{singletrans:00}) and (\ref{proofz:00}) taking $\mathbb{E}\big [Z_t[k]\big ] = \theta_t$. This Bayesian model is much simpler than the de Finetti model in Example \ref{deFinetti} where the likelihood in a transition is
\begin{equation*}
\int_{[0,1]}
\prod_{i=0}^1\prod_{j=0}^1{n_{i\cdot}\choose n_{ij}}p_t(i,j;\theta)^{n_{ij}^t}
\ \nu(d\theta).
\end{equation*}
Maximum likelihood can now be used to obtain an estimate $\widehat{\theta}_t$ in each transition. If $N$ is large this estimate will be approximately normally distributed with a variance proportional to $N^{-1}$.  Information about $\nu$ is then obtained by the estimates $\widehat{\theta_t}$ over transitions. 
There are two independent equations in the following in each transition when $\varphi\ne \gamma$. Taking $(i,j)=(0,0),(1,0)$ two equations from (\ref{like:343}) are
\begin{equation*}
\begin{aligned}
 p_0(\theta_t)
&:=
{\varphib}
\Big \{1 + \frac{{\alpha}}{{\alphab}}\big (1 - \frac{{\theta_t}}{{\alpha}} \big )
\Big \}=
\begin{cases}
1 - \theta_t& \varphi\leq \gamma\\
\frac{\varphib}{\gammab}(1-\theta_t) & \varphi > \gamma
\end{cases}
\\
p_1(\theta_t)
&:=
{\varphib}
\Big \{1 - \frac{{\alpha}}{{\alphab}}\big ( 1 - \frac{{\theta_t}}{{\alpha}} \big )
\frac{{\gammab}}{{\gamma}}\Big )
\Big \}=
\begin{cases}
1 - \varphi/\gamma + \theta_t\gammab/\gamma & \varphi\leq \gamma\\
\frac{\varphib}{\gamma}\theta_t & \varphi > \gamma.
\end{cases}
\end{aligned}
\label{like:101}
\end{equation*}
An estimate $\widehat{\theta}_t\in (0,1)$ is found by solving the derivative of the log likelihood equal to zero. That is solving 
\begin{equation}
\begin{aligned}
0=c_0(\varphi,\gamma)\Big (-n^t_{00}\frac{1}{p_0(\theta_t)} + n^t_{01} \frac{1}{1-p_0(\theta_t)}\Big )
+ c_1(\varphi,\gamma)\Big (n^t_{10}\frac{1}{p_1(\theta_t)}  -n^t_{11}\frac{1}{1-p_1(\theta_t)}\Big )
\end{aligned}
\label{like:106}
\end{equation}
where

\begin{equation*}
c_0(\varphi,\gamma) =
\begin{cases}
1&\varphi \leq \gamma\\
\frac{\varphib}{\gammab}&\varphi > \gamma
\end{cases},\quad
c_1(\varphi,\gamma) =
\begin{cases}
\frac{\gammab}{\gamma}&\varphi \leq \gamma\\
\frac{\varphib}{\gamma}&\varphi > \gamma
\end{cases}.
\end{equation*}
A small simulation study was made with $N=1000$ and $T=5000$.
\medskip

\vbox{
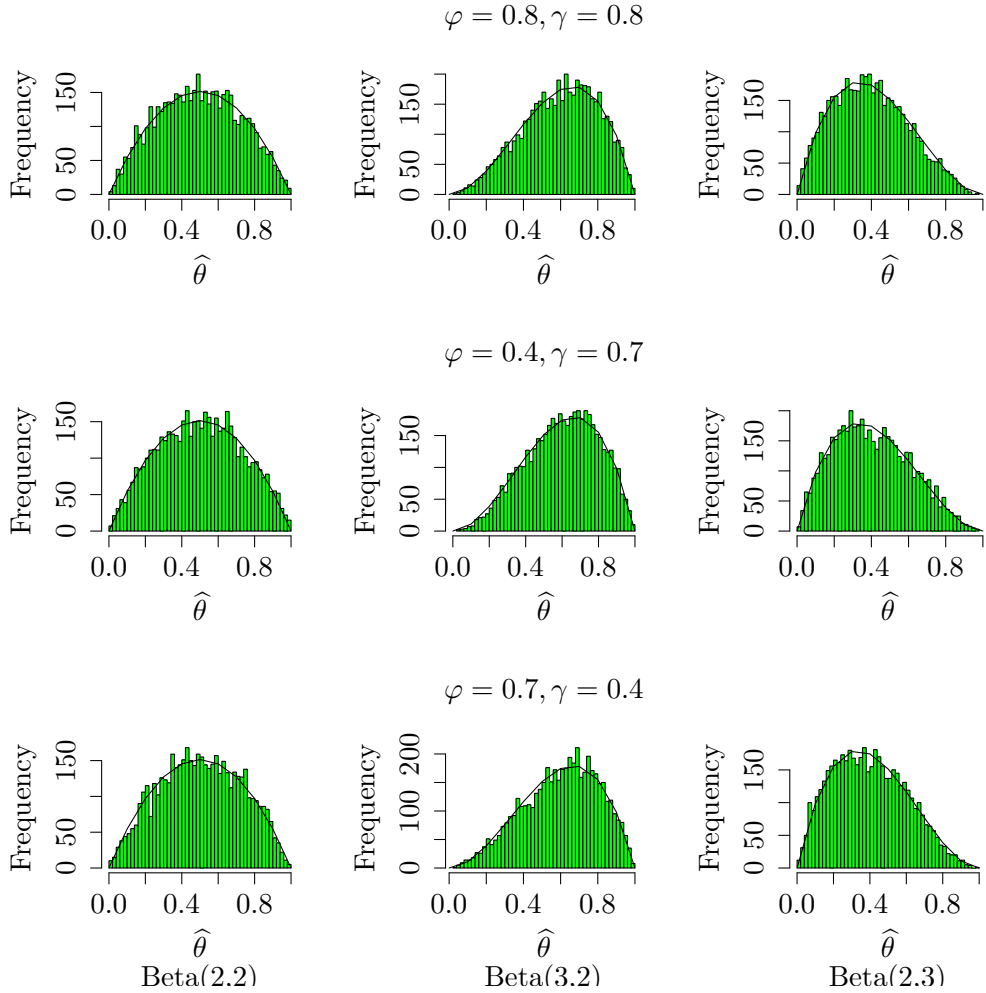
\emph{Histograms of estimated $\theta$ distributions, compared to true distribution curves.}\\
\noindent
\input{pic.tex}%[scale=.3]
}

\noindent
Parameter values $\varphi,\gamma$ shown are the same within rows. $\theta$ is chosen from Beta distributions with different parameters in columns. In each experiment the length of $X_t$ is $N=1000$ and $T=5000$ transitions are simulated with an initial distribution simulated from (\ref{Intro:2}). $\widehat{\varphi},\widehat{\gamma}$ estimates were very accurate and identical to true values to two significant figures. Histograms are from $\widehat{\theta}_t$ values estimated in each transition from (\ref{like:106}). They are compared with the expected curve from  Beta distributions. There is a close agreement between the histograms, which are estimated from sample path data, and exact theoretical curves.
\section*{Acknowledgement}\ {Andrea Collevecchio\rq{}s  work is partially supported by Australian Research Council grant  DP180100613 and  by Australian Research Council Centre of Excellence for Mathematical and Statistical Frontiers (ACEMS) CE140100049.}

\end{document}